\documentclass[11pt]{article}
\usepackage[utf8]{inputenc}
\usepackage{amsmath}
\usepackage{amssymb}
\usepackage{amsthm}
\usepackage{graphicx}
\usepackage{hyperref}
\usepackage{stmaryrd}
\usepackage{multirow}
\usepackage{rotate}
\usepackage{microtype}
\usepackage{courier}
\usepackage{color}
\usepackage[numbers,sort,compress]{natbib}
\usepackage[all]{hypcap}
\usepackage{doi}
\usepackage{enumitem}
\usepackage{tikz}
\usetikzlibrary{arrows}
\usetikzlibrary{positioning,calc}
\usepackage{tikzinclude}
\usepackage[ruled,vlined]{algorithm2e}
\usepackage{thmtools}
\declaretheoremstyle[headfont=\itshape, bodyfont=\normalfont, qed=$\triangle$, postheadspace=1em]{assumpstyle}
\declaretheorem[style=assumpstyle]{assumption}
\usepackage{mathtools}
\usepackage{geometry}
\theoremstyle{plain}
\newtheorem{lemma}{Lemma}
\newtheorem{proposition}{Proposition}
\newtheorem{theorem}{Theorem}
\newtheorem{corollary}{Corollary}
\newtheorem{remark}{Remark}

\hypersetup{
   unicode=true,
   pdftoolbar=false,
   pdfmenubar=false,
   pdffitwindow=false,
   pdfnewwindow=true,
   pdfkeywords={},
   colorlinks=true,
   linkcolor=black,
   citecolor=black,
   filecolor=black,
   urlcolor=black
}

\newcommand{\R}{\mathbb{R}}

\newcommand{\X}{\mathbf{X}}

\newcommand{\vetx}{\mathbf{x}}
\newcommand{\vety}{\mathbf{y}}
\newcommand{\vetz}{\mathbf{z}}

\newcommand{\vete}{\mathbf{e}}
\newcommand{\vetg}{\mathbf{g}}

\newcommand{\bb}{\begin{equation}}
\newcommand{\ee}{\end{equation}}
\newcommand{\ds}{\displaystyle}
\newcommand{\II}{\mathcal{I}}

\newcommand{\PX}{\mathcal{P}_{\mathbf{X}}}
\newcommand{\PXast}{\mathcal{P}_{\mathbf{X^{\ast}}}}

\newcommand{\gks}{\mathbf{g}_{s_{\ell}(k+\delta)}^{k+\delta-1}}
\newcommand{\gksq}{\mathbf{g}_{s_{\ell}(k+j(q))}^{k+ \delta-q}}
\newcommand{\eks}{\boldsymbol \epsilon_{s_{\ell}(k+ \delta)}^{k+\delta}}
\newcommand{\eksq}{\boldsymbol \epsilon_{s_{\ell}(k+j(q))}^{k+j(q)}}
\newcommand{\eksdq}{\boldsymbol \epsilon_{s_{\ell}(\delta k + j(q))}^{\delta k + j(q)}}

\newcommand\ceil[1]{\lceil#1\rceil}
\newcommand{\dxy}{\| \mathbf{x}^k - \mathbf{y} \|^2}

\newcommand{\EE}{\mathbb{E}}

\newcommand{\Newstar}{s_\ell(\delta n(k) ) \in \mathcal{R}_v}

\newcommand{\eliasc}[2]{{\color{black}{#2}}}
\newcommand{\rafael}{}
\newcommand{\rafaelb}{}
\newcommand{\rafaelc}{}
\newcommand{\eduardo}{}
\newcommand{\elias}{}
\newcommand{\eliasb}{}

\definecolor{darkblue}{rgb}{0,0,.65}

\title{A Markovian Incremental Stochastic Subgradient Algorithm}

\author{
Rafael Massambone\thanks{State University of Northern Paran\'{a} -- Centro de Ci\^{e}ncias Tecnol\'{o}gicas, C.P. 261,  86360--000, Bandeirantes, PR, Brazil. email: massambone@uenp.edu.br}
\and
Eduardo F. Costa\thanks{University of S\~{a}o Paulo -- Instituto de Ci\^{e}ncias Matem\'{a}ticas e de Computa\c{c}\~{a}o, C.P. 668, 13560--970, S\~{a}o Carlos, SP, Brazil. emails: \{efcosta,elias\}@icmc.usp.br}
\and
Elias S. Helou\footnotemark[2]
}

\begin{document}

\maketitle

\begin{abstract}
A stochastic incremental subgradient algorithm 
for the minimization of a sum of convex functions is
introduced. The method sequentially uses partial subgradient
information and the sequence of partial subgradients is
determined by a general Markov chain. This makes it suitable
to be used in networks where the path of information flow is
stochastically selected. We prove convergence of the
algorithm to a weighted objective function where the weights
are given by the \rafaelc{Cesàro} limiting probability distribution 
of the Markov chain. Unlike previous works in the literature, the Cesàro
limiting distribution is general (not necessarily uniform),
allowing for general weighted objective functions and
flexibility in the method.
\linebreak
\textbf{keywords:} Optimization algorithms, incremental subgradient method, stochastic optimization, Markov chains.
\end{abstract}

\section{Introduction}
Incremental subgradient algorithms~\cite{nedic01,helou09} can be efficient tools for the minimization of a sum of convex functions $f = \sum_{i = 1}^mf_i$. These methods use only the subgradient of a summand $f_i$ of the objective function at a time when updating the decision variables and, therefore, need to somehow move along all summands in order to converge to a minimizer of $f$. Several schemes have been used for the choice of the sequence of updates, starting with the most intuitive cyclic scheme~\cite{kibardin1979decomposition} going through each summand exactly once before repeating any. The cyclic scheme can be used with random permutations following each full cycle and this can be beneficial in practice. Another possibility is to choose uniformly and randomly the summand to be considered, which also seems to favor convergence speed in practice when compared to a deterministic cyclic scheme.

A common application of incremental subgradients methods is in optimization in networks. In this case each summand is processed by a network agent and, in order to avoid extra communication and processing costs and to take the network topology into consideration, it is usually necessary to impose restrictions on which summands can be selected for processing after some $f_i$ is processed. Usually the sequence must follow the neighborhood structure of the network, that is, summand $f_i$ can be followed by summand $f_j$ if and only if node $i$ is connected to node $j$. Under these circumstances, the randomizations of the processing order are not practically implementable without causing an unnecessary load in the network and, instead, a Markov chain can be used in order to select the next node from all neighbors of the current one~\cite{sundhar2009incremental, johansson2010}. \eliasb{Another useful application of incremental methods is in computed tomography~\cite{helou05,helou09,herman2009fundamentals} where the combination of fast initial convergence and computationally simple iterations is favoured over several digits accuracy.}

These Markovian methods have so far been proven to converge under the hypothesis of uniformity of the limiting probability distribution of the Markov chain. In the present paper we remove this restriction and prove convergence to the minimizer of an averaged objective function $\sum_{i = 1}^m\omega_i f_i$ where each $\omega_i$ correspond to the Cesàro limiting probability of node $i$ in the Markov chain, as intuitively should happen since $f_i$ will be processed with relative frequency of $\omega_i$.

\eduardo{Our setup allows e.g. to consider a Markov chain for which $\omega_j=0$ for some indexes $j$, so that the corresponding functions $f_j$ play a role in the initialization and initial steps (to speed up the method) while they do not interfere with the final solution. It is worth a mention that the Markov chain is usually constructed/designed by the user of the method, so both $\omega_i$ and the transition probabilities can be chosen at the user convenience.}
\elias{Not only we provide convergence results for the case of general limiting probability distributions, but our method is more general than other Markovian incremental algorithms because it allows parallel processing of certain subsets of summands of the objective function.} \rafael{In this sense, the proposed algorithm has similarities with the class of string-averaging incremental subgradient methods~\cite{helou14, censor2014string, oliveira16, oliveira19} for constrained convex optimization problems.} \elias{Indeed, the proposed method seems to be the first instance of Markovian string-averaging technique in the literature.} \eliasb{Interestingly, here again, computed tomography is an application for which string-averaging is useful~\cite{helou14,oliveira16,oliveira19}.}

\eduardo{Such a general setup brings some challenges for the convergence analysis,
which demanded for the development of Lemmas \ref{lem-eval-n1-2} and \ref{lemma-main-estimate} that are unparalleled in literature. We give more details
in Remark \ref{rem-comparing-with-sundhar2009incremental}.
Still, the main convergence result given in Theorem \ref{teorema1} is quite conventional and does not impose
restrictive / working assumptions. For example, Theorem \ref{teorema1} recovers  Theorem 4.3 of \cite{sundhar2009incremental}
in the context considered in that paper.}

The contributions of the present paper are, therefore, threefold: (a) it is shown that Markovian incremental stochastic subgradient algorithms converge under more general conditions on the Markov chain and on the stochastic gradient error than previously established in the literature; (b) the Markovian incremental stochastic subgradient algorithms are generalized to work with many parallel Markov chains in a string-averaging scheme, and (c) theoretical analysis of global convergence under mild assumptions is presented. We analyse both diminishing and constant stepsizes. \rafaelc{Our approach unifies different algorithms found in literature, such as the incremental (cyclic) subgradient method, Markov randomized incremental subgradient method and incremental (randomized) subgradient method \cite{nedic01,sundhar2009incremental,johansson2010}}. Moreover, we also provide experimental evidence of the usefulness of the method, although with no intention of pursuing benchmarking purposes.

The paper is organized as follows:
\eliasb{Subsection \ref{subsec:problem} describes precisely the problem we are aiming to solve.} Subsection \ref{sec.2} provides a complete description of the algorithm we propose;
Section \ref{sec.3} \eliasb{is dedicated to the theoretical analysis};
in Section \ref{sec.4} we examine a numerical example;
final considerations are given in Section \ref{sec.5}.

\section{Proposed Algorithm}

\subsection{Problem formulation}\label{subsec:problem}

We consider a network of $m$ agents that are indexed by $\{1,\dots, m\}$, and we denote $\II := \{1,\dots, m\}$. The network has a
static topology that is given by the directed graph $G := (\II, \mathcal{E})$, where $\mathcal{E} \subset \II^2$ is the set of links in the network.
We have $(i, j) \in \mathcal{E}$ if agents $i$ and $j$ can communicate with each other.
In this paper, the network goal is to solve the following optimization problem:
\begin{equation} \label{main-problem} \begin{array}{c} \ds \vetx \in \arg \min f(\vetx)
\\
\mbox{s.t.} \quad \vetx \in \X \subseteq \mathbb{R}^n, \end{array} \end{equation}
where
\begin{description}[listparindent=1.5em, leftmargin=0cm]
\item[(i)] $ \ds f(\vetx) := \sum_{i \in \II} \omega_i f_{i}(\vetx), \quad \omega_i \ge 0, \quad \sum_{i \in \II} \omega_i = 1;$
\item[(ii)] $f_i : \R^n \rightarrow \R$ are convex for all $i \in \II$.
\end{description}

\rafael{The literature on methods to solve the optimization problem (\ref{main-problem}) is very wide, so that not only incremental methods have been studied for this purpose. We can highlight, for example, distributed computation techniques \cite{nedic09distributedIEEETAC, sundhar2010distributed, nedic10IEEETAC, jakovetic14IEEETAC, nedic15IEEETAC} where each agent $i$ minimizes its own objective $f_i$ by exchanging information with other agents on the network. Distributed algorithms and optimization models have been studied for a long time, mainly due to potential in applications such as sensor networks \cite{rabbat04} and distributed control \cite{sundhar2009distributedIEEE}.} \eliasb{One of the main disadvantages of distributed methods compared to the \emph{Markovian Incremental Stochastic Subgradient Algorithm} (MISSA) that we propose here is that they require more sophisticated oracles (the minimizer instead of a stochastic subgradient) which might not always be available or be practical.}

\subsection{Algorithm description} \label{sec.2}
\begin{figure}[t!]
\centering
\begin{tikzpicture}[->,>=stealth',shorten >=1pt,auto,node distance=1.5cm,
thick, main node/.style={transform shape, circle, minimum size=3mm, fill=blue!10, draw}]
\node (S1) at (0.7, 2) {$\mathcal{R}_1$};
\node[main node] (i11) at (3, 3) {$1$};
\node[main node] (i12) [below left of=i11] {$2$};
\node[main node] (i13) [below right of=i12] {$3$};
\node[main node] (i14) [below right of=i11] {$4$};

\node[main node] (i21) [right= 3cm of i11] {$8$};

\node[main node] (i31) [below= 1.5cm of i13] {$9$};

\node (S4) at (8.8,-2) {$\mathcal{R}_2$};
\node[main node] (i41) [right= 2.15cm of i31] {$5$};
\node[main node] (i44) [below right of=i41] {$7$};
\node[main node] (i42) [above right of=i44] {$6$};

\path[every node/.style={font=\sffamily\small}]
(i14) edge [right] node[left] {} (i11)
(i11) edge [right] node[left] {} (i13)
(i11) edge [bend left] node[left] {} (i14)
(i12) edge [bend right] node[left] {} (i13)
(i13) edge [bend left] node[left] {} (i11)
(i13) edge [right] node[left] {} (i12)
(i12) edge [right] node[left] {} (i14);

\path[every node/.style={font=\sffamily\small}]
(i41) edge [right] node[left] {} (i44)
(i42) edge [right] node[left] {} (i41)
(i44) edge [right] node[left] {} (i42);

\path[every node/.style={font=\sffamily\small}]
(i21) edge [loop] node[left] {} (i21)
(i21) edge [bend left] node[left] {} (i13)
(i21) edge [right] node[left] {} (i42);

\path[every node/.style={font=\sffamily\small}]
(i31) edge [bend right] node[left] {} (i11)
(i31) edge [right] node[left] {} (i44);
\end{tikzpicture}
\caption[A network example.]{A disconnected graph as an example of a network topology containing two recurrent class and two transient states.} \label{network_example}
\end{figure}
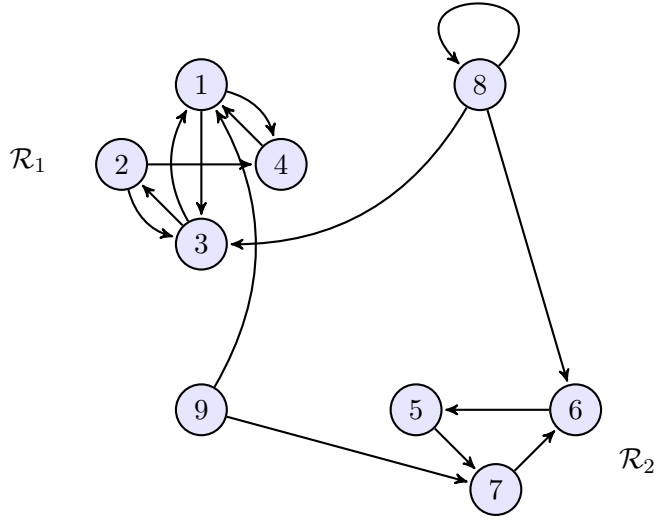

The algorithm we propose in this paper, named as MISSA, is presented in this section.
Consider $M$ Markov chains indexed by $\{1,\ldots,M\}$ and we
denote the index set by $\mathcal M=\{1,\ldots,M\}$, so that for each $\ell\in \mathcal M$
we have that $\{s_\ell(k),k\geq 0\}$ is a time-homogeneous Markov chain with
state space $\mathcal I$ and transition probabilities $\text{Prob}(s_\ell(k+1)=j|
s_\ell(k)=i)=[P]_{ij}$, where $P$ is called the transition probability matrix and is common to all $M$
Markov chains\footnote{\eduardo{$P$ can be seen as a
parameter of the method, to be selected by the user aiming at a desired Cesàro limiting
distribution and/or transient states for accelerating the algorithm.}}.
According to \cite[Theorem 3.14]{cinlar2013introduction}, a simple re-ordering of the Markov states allow to
obtain the following form for $P$,
\begin{equation} \label{trans_prob_matrix}\small
P = \left[ \begin{array}{ccccc}
P_1 &     &        &  & \\
    & P_2 &        &  & \\
    &     & \ddots &  & \\
    &     &        & P_N & \\
Q_1 & Q_2 & \cdots & Q_N & T
\end{array} \right],
\end{equation}
where for each $v=1,\dots,N$, the matrix  $P_v$ is of
dimension $m(v) \times m(v)$, $Q_v$ is a $u \times m(v)$ matrix whereas $T$ is a $u \times u$ matrix.
The above structure of $P$ makes evident which states are recurrent
\cite{cinlar2013introduction}: the rows of $P$ corresponding to $P_v$ for a specific $v$
form a set of Markov states
that are all recurrent and of same period \cite{cinlar2013introduction};
we denote the period by $\delta_{v}$
and this set of Markov states by $\mathcal{R}_v$.
The rows of $P$ corresponding to the matrix $T$ form a set of transient Markov states,
which we denote by $\mathcal{T}$.  Note that $s_\ell(k)$ takes values on $\II = (\bigcup_{v=1}^{N} \mathcal{R}_{v}) \cup \mathcal{T}$.
As an example,
consider a network composed of $9$ agents, so that $\II := \{1,\dots, 9\}$,
and $4$ Markov chains $s_1(k), \dots, s_4(k)$, whose probability matrices are given by a common matrix
$$\small P=\begin{bmatrix}
0 & 0 & 0.5 & 0.5 &&&&&&\\
0 & 0 & 0.3 & 0.7 &&&&&&\\
0.2 & 0.8 & 0 & 0  &&&&&&\\
 1 &  0 & 0 & 0  &&&&&& \\
&&&&0&0&1 &&&\\
&&&&1&0&0 &&&\\
&&&&0&1&0 &&&\\
0&0&0.1&0&0&0.2&0&0.7&0\\
0.1&0&0&0&0&0&0.9&0&0
\end{bmatrix}.
$$
In this example we have two \eduardo{irreducible} sets of recurrent states $\mathcal{R}_1 = \{1,2,3,4\}$, $\mathcal{R}_2 = \{5,6,7\}$ and a set of transient states
$\mathcal{T} = \{8,9\}$.
States in $\mathcal{R}_1$ are periodic with period $\delta_1=2$, while states in $\mathcal{R}_2$ are of period $\delta_2 = 3$.
It is worth mentioning that $P^\delta$ is the transition probability matrix of an
aperiodic Markov chain, where $\delta \,= \, \text{LCM}(\delta_1,\ldots, \delta_N)$
(where LCM means ``least commom multiple'');
in the example, $\delta \,= \, \text{LCM}(2,3)=6$ so that $P^6$ gives the transition
probabilities of an aperiodic chain, or equivalently
we can say that the Markov chains $\{s_1(6k),k\geq 0\},\ldots,\{s_4(6k),k\geq 0\}$ are aperiodic. Figure \ref{network_example} shows the network topology $G$ of this example.

\rafael{Let us define \elias{$\boldsymbol \pi_{\ell}^{0} \in \mathbb R^m$} as the initial \elias{probability} distribution for $s_{\ell}(k)$, for all $\ell \in \mathcal{M}$}\elias{, i.e., $\text{Prob}( s_\ell( 0 ) = i ) = [\boldsymbol\pi_{\ell}^{0}]_i$}.
We assume that every agent is relevant to the algorithm in the sense that it
is visited at \eliasc{some iteration}{infinitely many iterations} $k$ with
positive probability (otherwise the agent could be removed from network with no loss).
Returning to the example above, this prevent us from setting\elias{, for example, }$\boldsymbol \pi_{1}^{0} = \boldsymbol \pi_{2}^{0} = \boldsymbol \pi_{3}^{0}
= \boldsymbol \pi_{4}^{0} = [0.5, 0.5, 0, \dots, 0]^{T}$
because the states in $\mathcal{R}_2 = \{5,6,7\}$ would never be visited by any Markov chain. 
An example of valid initial distributions is
$\boldsymbol \pi_{1}^{0} = [0.5, 0.5, 0, \dots, 0]^{T}$,
$\boldsymbol \pi_{2}^{0} = [0, 0, 0, 0, 1, 0, 0, 0, 0]^{T}$,
$\boldsymbol \pi_{3}^{0} = [0, 0, \dots, 1, 0]^{T}$
and $\boldsymbol \pi_{4}^{0} = [0, 0, \dots, 0, 1]^{T}$.

Some strong hypothesis, like the distribution converging to $[1/m, \dots, 1/m]^{T}$, are void in our setup - in fact, we do not require existence nor uniqueness of a limiting
distribution for $\{s_\ell(k)\}$, $\ell \in \mathcal{M}$.
We resort to the limiting distribution in the
sense of Ces{\`a}ro \cite{spreen1981}:
\begin{equation} \label{cesaro1}
\bar P :=  \lim_{k \to \infty} \frac{1}{k} (P^{0} + P^{1} + \dots + P^{k-1}).
\end{equation}
If we define $M_k = (P^k+\cdots+ P^{k+\delta-1})$ then the definition given in
\eqref{cesaro1} is equivalent to
$\bar P =  \frac{1}{\delta} \lim_{k \to \infty} \frac{1}{k} (M_{\delta 0} + \cdots + M_{\delta (k-1)})
= \frac{1}{\delta} \lim_{k \to \infty} \frac{1}{k} M_0 (P^0 + P^\delta + \cdots + P^{\delta (k-1)}).$
Since $P^{\delta}$ is aperiodic \cite[Theorem 3.7, p.162]{cinlar2013introduction}, then
$P^{\delta k}$ converges, yielding
$\bar P = \frac{1}{\delta} \lim_{k \to \infty} M_0 P^{\delta k},$ and this in turn leads to
\begin{equation} \label{cesaro2}
\bar P = P_{\delta} := \frac{1}{\delta} \, \lim_{k \to \infty} (P^{k} + P^{k+1} + \dots + P^{k + \delta - 1}).
\end{equation}
Now we can use equation \eqref{cesaro2} to obtain for all $\ell \in \mathcal{M}$ and any initial distribution $\boldsymbol \pi_{\ell}^{0}$, a Ces{\`a}ro limit distribution
$\boldsymbol \pi_{\ell}^{\infty}$ as
\begin{equation} \label{limiting_distribution}
(\boldsymbol \pi_{\ell}^{\infty})^T = (\boldsymbol \pi_{\ell}^{0})^{T} P_{\delta}.
\end{equation}
\rafaelc{Additionally, we define $\Delta := \lim_{k\rightarrow\infty} P^{\delta k}$ and $\Phi( k, t ) := P^{k-t}$ with $k>0$, $t \geq 0$ and $k > t$.}

\elias{Summarizing} the method we propose,
at iteration $k$, we compute $\ell$ independent subiterations by using a stochastic subgradient of $f_{s_\ell(k)}$ and then we perform an averaging of these subiterations. \eliasb{It is useful to remark that although averaging is a possibility encompassed by our method, it is not mandatory. Indeed, if $|\mathcal M| = 1$ we recover the fully incremental method and still in this case our convergence results are novel because of the generality of the allowed Markov chain.}
The method \eliasc{MISSA}{Markovian Incremental Stochastic Subgradient Algorithm (MISSA)} is detailed in Algorithm \ref{missa}.
\begin{algorithm}[h!] \caption[MISSA]{\textsc{Markovian Incremental Stochastic Subgradient Algorithm - MISSA}} \label{missa}
\SetKwData{Left}{left}\SetKwData{This}{this}\SetKwData{Up}{up}
\SetKwFunction{Union}{Union}\SetKwFunction{FindCompress}{FindCompress}
\SetKwInOut{Input}{input}\SetKwInOut{Output}{output}
\Input{An integer $K > 0$;
an initial vector $\vetx^0 \in \X$; a set $\{s_\ell(k)\}$ of Markov chains (with transition probability matrix $P$ given by (\ref{trans_prob_matrix}));
a sequence $\{ \lambda_k \}$ of positive stepsizes.}
\For{$k=0$ \KwTo $K$}{
For each $\ell \in \mathcal{M}$, compute
\begin{align}
\label{missa1} \ds \vetx_{\ell}^{k+1} &= \vetx^k - \lambda_{k}(\vetg_{s_\ell(k+1)}^{k} + \boldsymbol \epsilon_{s_\ell(k+1)}^{k+1}),
\end{align}
where $\vetg_{s_\ell(k+1)}^{k} \in \partial f_{s_\ell(k+1)}(\vetx^k)$ and $\boldsymbol \epsilon_{s_\ell(k+1)}^{k+1}$ is a random noise vector.

Use the subiterations $\vetx_{\ell}^{k+1}$ obtained in (\ref{missa1}) to compute:
\begin{align}
\label{missa2} \ds \vetx^{k+1} &= \PX\left( \frac{1}{|\mathcal{M}|} \sum_{\ell \in \mathcal{M}} \vetx_{\ell}^{k+1}\right),
\end{align}
where $\PX$ is the Euclidean projection onto $\X$.
}
\Return{$\vetx^{K+1}$.}
\end{algorithm}

\section{Convergence analysis} \label{sec.3}

\subsection{Assumptions} \label{subsec:assumptions}

\rafaelc{\emph{Notations}:} \eliasc{}{Let $\vetx \in \mathbb R^n$, then we will use the notations $\| \vetx \| := \sqrt{\sum_{i = 1}^n|x_i|^2}$ and $\| \vetx \|_\infty := \max_{i \in \{ 1, 2, \dots, n \}}| x_i |$. \rafaelc{For any $x \in \mathbb{R}$, $\lceil x \rceil$ is the smallest integer larger than or equal $x$}. For matrices $A\in \mathbb R^{m \times n}$, we use the induced norms $\| A \| := \sup_{\| \vetx \| = 1}\| A\vetx \|$ and $\| A \|_\infty := \sup_{\| \vetx \|_\infty = 1}\| A\vetx \|_\infty$.} We analyze the Algorithm \ref{missa} under the following assumptions.

\begin{assumption} \label{assump1} $\X$ is nonempty, convex and compact.
\end{assumption}
Assumption \ref{assump1} and convexity of each $f_i$ over $\R^n$ imply that the subgradients of $f_i$ are bounded over $\X$ for each $i \in \II$, i.e., there is $C > 0$ such that for all $\vetx \in \X$ and $i \in \II$:
\begin{equation} \label{subgrad_boundedness} \|\vetg_i \| \leq C, \quad \vetg_i \in \partial f_i(\vetx). \end{equation}

The stochastic error when computing a subgradient $\vetg_{i}^{k}$ is denoted by $\boldsymbol \epsilon_{i}^{k+1}$.
We define $\mathcal{G}_k$ as the $\sigma$-algebra generated by $\vetx^0$, the subgradient errors $\{ \boldsymbol \epsilon_{1}^{1}, \dots, \boldsymbol \epsilon_{m}^{k} \}$ and $\{ s_{\ell}(1), \dots, s_{\ell}(k) \}$ for all $\ell \in \mathcal{M}$.
In our analysis, we suppose that the second moments of the subgradient errors are bounded uniformly over the agents, conditionally on the available information until iteration $k$ (and represented by $\mathcal{G}_k$).

\begin{assumption} \label{assump2} There is a scalar sequence $\left\{\nu_k\right\}$ such that
$\EE [ \| \boldsymbol \epsilon_{i}^{k+1} \|^2 | \,  \mathcal{G}_k ] \leq \nu_{k+1}^{2}$ for all $i \in \II$ and $k\geq 0$ with probability $1$.
\end{assumption}

The main convergence result requires additional hypotheses
on the bounds of the subgradient errors, which can be seen
as a natural extension of the ones in
\cite[Theorem 4.3]{sundhar2009incremental} to the case of periodic Markov chains.
\begin{assumption} \label{assump3}
Consider the scalar sequence $\left\{\nu_k\right\}$ given in Assumption \ref{assump2}.
We assume that
\begin{equation}
\begin{aligned}
\nu := \sup_{k \geq 0} \nu_k < \infty,
\end{aligned}
\end{equation}
and for any $2/3 < \xi \leq 1$, we have
\begin{equation} \label{assump_stepsize-moment}
\begin{aligned}
\sum_{k=0}^{\infty} \frac{\nu_{\delta k + \delta - q + 1}}{(k+1)^\xi} < \infty, \,\, \forall \, q=1,\dots,\delta.
\end{aligned}
\end{equation}
\end{assumption}

\subsection{Analysis}

In \eliasb{the remainder of} this section, we analyze the convergence of MISSA in order to solve problem (\ref{main-problem}) with a specific choice for scalars $\omega_i$. We define for each $i \in \II$,
\begin{align} \label{weights}
\omega_i = \frac{1}{|\mathcal{M}|} \sum_{\ell \in \mathcal{M}} [\boldsymbol \pi_{\ell}^{\infty} ]_i,
\end{align}
where $\boldsymbol \pi_{\ell}^{\infty}$ is the Cesàro limit distribution given by equation (\ref{limiting_distribution}).
This means that we are establishing larger weights for the local functions $f_i$ associated to the agents $i$ with higher
probability of use for long iterations $k$, that is, when $k \to \infty$. As a consequence, we have that the $u$ transient states
do not affect the objective function $f$, because $\omega_j = 0$ for all $j=1,\dots,u$. \rafaelc{Our main results are given in the following theorems. Define}
\begin{equation*}f^{\ast} = \inf_{\vetx \in \X} f(\vetx) \quad \mbox{and} \quad \X^{\ast} = \{\vetx \in \X \, | \, f(\vetx) = f^{\ast} \}.\end{equation*}

\begin{theorem} \label{teorema1} Suppose that Assumptions \emph{\ref{assump1}}--\emph{\ref{assump3}} hold.
Consider that for any $t \in \mathbb{N}$,
\begin{equation} \label{stepsize-theorem}
\lambda_{\delta t} = \lambda_{\delta t + 1} = \dots = \lambda_{\delta (t+1) - 1} = \frac{a}{(t+1)^\xi},
\end{equation}
where $a > 0$ and $2/3 < \xi \leq 1$.
Then, if $\{ \vetx^k \}$ is generated by Algorithm \emph{\ref{missa}} we have, with probability $1$,
\begin{equation*}\ds \lim_{k \to \infty} \inf f(\vetx^{\delta \eliasc{n(k)}{k}}) = f^{\ast} \quad \mbox{\textit{and}} \quad \lim_{k \to \infty} \inf d(\vetx^{\delta \eliasc{n(k)}{k}}, \, \X^{\ast})=0,\end{equation*}
where $d(\vetx, \, \X^{\ast}) := \|\vetx - \PXast(\vetx) \|$.
\end{theorem}

\eliasc{}{%
\begin{theorem} \label{teorema2}
   Suppose that Assumptions \emph{\ref{assump1}}~and~\emph{\ref{assump2}} hold.
   Consider that for any $k \in \mathbb{N}$,
   \begin{equation*}
   \lambda_k \equiv \lambda > 0, \quad\text{and}\quad \nu \ge \nu_k,
   \end{equation*}
   thus, if $\{ \vetx^k \}$ is generated by Algorithm \emph{\ref{missa}} then, for small $\lambda$,
   \begin{equation*}
      \liminf_{k \to \infty}\EE[ f(\vetx^{\delta k}) - f^{\ast}] \le \frac{\delta( \nu + C )}\theta\lambda\ln\frac1\lambda + A\lambda + M\nu,
   \end{equation*}
   for some \rafaelc{positive scalars $A,M$ and $\theta$}. For large $\lambda$ we have
   \begin{equation*}
      \liminf_{k \to \infty}\EE[ f(\vetx^{\delta k}) - f^{\ast}] \le \lambda\left[\frac{\delta( C + \nu )}\theta + \frac{( C + \nu )^2}2\right] \varrho e^{-\theta} + M\nu.
   \end{equation*}
\end{theorem}%
}

Before we proceed to \eliasc{\eliasb{prove this result}}{prove these theorems}, we need two important auxiliary results. \eliasb{The first of these results shows that the averaged expected value over $\mathcal R_v$ of any parcel $f_i$, $i \in \mathcal R_v$ of the objective function at a given iteration becomes less dependent of the next state of the Markov chain as iterations proceed.} \rafael{For that, we define the vector $\boldsymbol \pi_{\delta,v}$ as any row of $P_\delta$ such that the index of the row corresponds to some state in $\mathcal{R}_v$. Notice that, by the definition of $P_\delta$ in \eqref{cesaro2}, all these rows are equal. With such definition, for any $v=1,\dots,N$ and $h \in \mathcal{R}_v$,
\begin{equation}\label{eq-def-linha-Pdelta}
\vete_h^T [P_\delta]_i = [\boldsymbol \pi_{\delta,v}]_i,\eliasc{}{\quad i \in \{1, 2, \dots, m\},}
\end{equation}
where $\vete_h$ is the vector with $1$ in the $h$ coordinate and zero in the others and $[P_\delta]_i$ denotes the $i$-th column of $P_\delta$.}
\rafael{
\begin{lemma} \label{lem-eval-n1-2}
Let \eliasc{$n(k)$}{k} be a non-negative integer\eliasc{}{ and consider another integer $n( k ) \in [ 0, k - 1 ]$}. Define for each $q=1,\dots,\delta$, $j(q) = \delta - q + 1$.
For a given $\ell \in \mathcal{M}$ and a given $v$, $1 \leq v\leq N$,
there are positive scalars $\eta$ and $\beta$ such that, for any $\vety \in \X$,
we have
\begin{equation*}
\begin{aligned}
&\sum_{q=1}^\delta \sum_{i,h \in \mathcal{R}_v}  \EE[ ( f_i(\vetx^{\delta n(k)}) - f_i(\vety) )  1_{\{s_\ell(\delta n(k) +1) = h,\Newstar\}} ] [ P^{\delta (k-n(k)) + j(q) - 1} ]_{h i} \geq
- \eta e^{-\beta (k-n(k))} \\
& \quad +
\delta  \sum_{i \in \mathcal R_v}
\EE \Big[ ( f_i(\vetx^{\delta n(k)}) - f_i(\vety) )  1_{\{\Newstar\}} \Big] [\boldsymbol \pi_{\delta,v}]_i.
\end{aligned}
\end{equation*}
\end{lemma}
\begin{proof}
For ease of notation in this proof we write (keeping fixed $k > 0$):
\begin{equation*}\begin{aligned}
{\bf P}_{\delta}^{k} &= P^{\delta( k - n(k) ) - 1},\\
\Theta_{\ell,v}^{k} & = \sum_{q=1}^\delta \sum_{i,h \in \mathcal{R}_v}  \EE[ ( f_i(\vetx^{\delta n(k)}) - f_i(\vety) ) 1_{\{s_\ell(\delta n(k) +1) = h,\Newstar\}} ] [ P^{\delta (k - n(k) ) + j(q) - 1} ]_{h i}.
\end{aligned}\end{equation*}
If we note that
$\sum_{q=1}^\delta [ P^{\delta (k-n(k)) + j(q)  - 1} ]_{h i}
=\sum_{q=1}^\delta \vete_h^T [{\bf P}_{\delta}^{k} P^{j(q)}]_i = \vete_h^T [{\bf P}_{\delta}^{k} \sum_{q=1}^\delta P^{j(q)}]_i$, we can write
\begin{equation*}
\begin{aligned}
\Theta_{\ell,v}^{k} &= \sum_{i,h \in \mathcal R_v}
\EE \Big[ ( f_i(\vetx^{\delta n(k)}) - f_i(\vety) ) 1_{\{s_\ell(\delta n(k)+1)=h,\Newstar\}} \Big] \vete_h^T [{\bf P}_{\delta}^{k} \sum_{q=1}^\delta P^{j(q)}]_i
\\ & = \sum_{i,h \in \mathcal R_v}
\EE \Big[ ( f_i(\vetx^{\delta n(k)}) - f_i(\vety) ) 1_{\{s_\ell(\delta n(k)+1)=h,\Newstar\}} \Big] \Big(\vete_h^T \Big[ {\bf P}_{\delta}^{k} \sum_{q=1}^\delta P^{j(q)} \Big]_i + \vete_h^T [\delta P_\delta]_i - \vete_h^T [\delta P_\delta]_i \Big)
\end{aligned}
\end{equation*}
thus leading to
\begin{equation*}
\begin{aligned}
\Theta_{\ell,v}^{k} &= \sum_{i,h \in \mathcal R_v}
\EE \Big[ ( f_i(\vetx^{\delta n(k)}) - f_i(\vety) ) 1_{\{s_\ell(\delta n(k))+1=h,\Newstar\}} \Big] \vete_h^T [\delta P_\delta]_i
\\ & \quad +
\sum_{i,h \in \mathcal R_v}
\EE \Big[ ( f_i(\vetx^{\delta n(k)}) - f_i(\vety) ) 1_{\{s_\ell(\delta n(k))+1=h,\Newstar\}} \Big] \big(\vete_h^T \Big[{\bf P}_{\delta}^{k} \sum_{q=1}^\delta P^{j(q)}\Big]_i - \vete_h^T [\delta P_\delta]_i\big).
\end{aligned}
\end{equation*}
Then,
\begin{equation*}
\begin{aligned}
\Theta_{\ell,v}^{k} &\geq  \sum_{i,h \in \mathcal R_v}
\EE \Big[ ( f_i(\vetx^{\delta n(k)}) - f_i(\vety) ) 1_{\{s_\ell(\delta n(k))+1=h,\Newstar\}} \Big] \vete_h^T [\delta P_\delta]_i
\\ &
\quad -\sum_{i,h \in \mathcal R_v}
\Big| \EE \Big[ ( f_i(\vetx^{\delta n(k)}) - f_i(\vety) ) 1_{\{s_\ell(\delta n(k))+1=h,\Newstar\}}  \Big] \Big| \Big| \vete_h^T \Big[{\bf P}_{\delta}^{k} \sum_{q=1}^\delta P^{j(q)}\Big]_i - \vete_h^T [\delta P_\delta]_i \Big|.
\end{aligned}
\end{equation*}
Notice that
\begin{align*}
\Big| \vete_h^T \Big[{\bf P}_{\delta}^{k} \sum_{q=1}^\delta P^{j(q)}\Big]_i - \vete_h^T [\delta P_\delta]_i \Big| & \leq \| \vete_h \|_{\infty} \Big\| {\bf P}_{\delta}^{k} \sum_{q=1}^\delta P^{j(q)} -  \delta P_\delta \Big\|_{\infty} \\
& = \Big\| {\bf P}_{\delta}^{k} \sum_{q=1}^\delta P^{j(q)} -  \Delta \sum_{q=1}^\delta P^{j(q)} \Big\|_{\infty}  \\
& = \Big\| \sum_{q=1}^\delta \Phi(\delta(k-n(k))-1 + j(q), \delta 0 + 0) - \Delta P^{j(q)} \Big\|_{\infty} \\
&\leq \sum_{q=1}^\delta \alpha e^{- \beta (k - n(k))}, \,\, (\alpha, \, \beta > 0),
\end{align*}
where $\alpha$ and $\beta$ \eliasc{the preceding inequality is}{are} due \eliasc{to Lemma~\ref{lema3} and}{} to Lemma~\ref{lema4}-\textbf{(ii)} \eliasb{(see Appendix)}.
Therefore, by using \eqref{eq-def-linha-Pdelta}, rearranging the terms and substituting above,
\begin{equation*}
\begin{aligned}
\Theta_{\ell,v}^{k} &\geq \delta [\boldsymbol \pi_{\delta,v}]_i \sum_{i,h \in \mathcal R_v} \EE \Big[ ( f_i(\vetx^{\delta n(k)}) - f_i(\vety) ) 1_{\{s_\ell(\delta n(k)+1)=h,\Newstar\}} \Big]
\\ & \quad - \sum_{i,h \in \mathcal R_v} \Big| \EE \Big[ ( f_i(\vetx^{\delta n(k)}) - f_i(\vety) ) 1_{\{s_\ell(\delta n(k)+1)=h,\Newstar\}}  \Big] \Big| \alpha \delta e^{-\beta (k-n(k))} \\
& \geq \delta [\boldsymbol \pi_{\delta,v}]_i \sum_{i\in \mathcal R_v} \EE \Big[ ( f_i(\vetx^{\delta n(k)}) - f_i(\vety) ) \sum_{h \in \mathcal R_v} 1_{\{s_\ell(\delta n(k)+1)=h,\Newstar\}} \Big] - \sum_{i,h \in \mathcal R_v} \bar M \alpha \delta e^{-\beta (k-n(k))},
\end{aligned}
\end{equation*}
where we have used the compactness of $\X$ and convexity of
$f_i$ to guarantee the existence of $\bar M > 0$ such that
$|f_i(\vetx^{\delta n(k)}) - f_i(\vety)|\leq \bar M$, yielding
\begin{equation*}\begin{aligned}
\Big| \EE \Big[ ( f_i(\vetx^{\delta n(k)}) - f_i(\vety) ) 1_{\{s_\ell(\delta n(k)+1)=h,\Newstar\}}  \Big] \Big| 
& \leq \EE \Big[\Big|  ( f_i(\vetx^{\delta n(k)}) - f_i(\vety) ) 1_{\{s_\ell(\delta n(k)+1)=h,\Newstar\}}  \Big|\Big] \\
& \leq  \bar M \EE \Big[1_{\{s_\ell(\delta n(k)+1)=h,\Newstar\}}\Big] \\
& \leq \bar M.
\end{aligned}
\end{equation*}
The fact that $\mathcal R_v$ is a closed set of Markov states yields
$$\sum_{h \in \mathcal R_v} 1_{\{s_\ell(\delta n(k)+1)=h,\Newstar\}}=1_{\{\Newstar\}},$$
hence
\begin{equation*}
\begin{aligned}
\Theta_{\ell,v}^{k} &\geq
\sum_{i\in \mathcal R_v}
\EE \Big[ ( f_i(\vetx^{\delta n(k)}) - f_i(\vety) )  1_{\{\Newstar\}} \Big]
\delta [\boldsymbol \pi_{\delta,v}]_i - \sum_{i,h \in \mathcal R_v} \bar M \alpha \delta e^{-\beta (k-n(k))}
\\ & \geq \sum_{i\in \mathcal R_v} \EE \Big[ ( f_i(\vetx^{\delta n(k)}) - f_i(\vety) )  1_{\{\Newstar\}} \Big]
\delta [\boldsymbol \pi_{\delta,v}]_i - \eta e^{-\beta (k-n(k))},
\end{aligned}
\end{equation*}
where $\eta = \eliasc{v}{m}^2 \bar M \alpha \delta$.
\end{proof}
}

Notice that we select two iterations, namely $k$ and a previous one $n( k )$, and then base our argument in comparing iteration $\vetx^{\delta k}$ with iteration $\vetx^{\delta n( k )}$. The reason why we do so is that we need to perform the analysis on later iterations, that is, we need that both $k \to \infty$ and $n( k ) \to \infty$ (such as when comparing iteration $k$ against iteration $k - 1$ in regular non-stochastic convergence analysis for subgradient methods), but this is not sufficient because the convergence in the markovian method also depends on analysing a long enough part of the chain, that is, the argument also requires $k - n( k ) \to \infty$ in a specific way, as seen in the proof of Theorem~1.

The next lemma provides an important estimate for convergence analysis of the Algorithm \ref{missa}. It gives means to establish an expected reduction in distance from the current point to any other given point with smaller objective function value as long as the stepsize is large enough, the stochastic error in subgradient computation is small and we have made sufficient passes over the Markov chain.

\begin{lemma} \label{lemma-main-estimate}
Suppose that Assumptions \emph{\ref{assump1}} and \emph{\ref{assump2}} hold. Consider that the stepsizes satisfy, for all $k \geq 0$
\begin{equation} \label{stepsize1}
\lambda_{\delta (k+1) - q} = \lambda_{\delta k}, \quad \forall \, q = 1, \dots, \delta.
\end{equation}
Let $n(k)$ be a nonnegative integer sequence satisfying $n(k) \leq k$ for all $k \geq 0$. If $\{\vetx^k \}$ is the sequence generated by \emph{Algorithm \ref{missa}} then, for all $\vety \in \X$, there are positive scalars $\eliasc{\beta}{\theta}$, $\varrho$ and $C$ such that
\begin{align} \label{eq-lemma-main-estimate}
\ds \EE[\| \vetx^{\delta (k+1)} - \vety \|^2]  &\leq \EE[ \| \vetx^{\delta k} - \vety \|^2] +  2\lambda_{\delta k} \sum_{q=1}^{\delta} \sum_{p=\delta n(k)}^{\delta (k+1)-q-1} \lambda_{p}(C + \nu_{p+1}) -  2 \delta \lambda_{\delta k} \EE[ f(\vetx^{\delta n(k)}) - f(\vety)] \nonumber \\
									& \quad + 2 \varrho \delta \lambda_{\delta k} (e^{- \eliasc{\beta}{\theta} n(k)} + e^{-\eliasc{\beta}{\theta} (k-n(k))}) +  2\lambda_{\delta k} \sum_{q=1}^{\delta} \nu_{\delta k + j(q)} \EE[ \| \vetx^{\delta (k+1)-q} - \vety \| ] \nonumber \\
									& \qquad + \lambda_{\delta k}^{2} \sum_{q=1}^{\delta}( \nu_{\delta k + j(q)} + C )^2.
\end{align}
\end{lemma}
\begin{proof}
Initially, due the non-expansivity of the projection $\PX$, equation (\ref{missa2}) and convexity of squared norm we have
\begin{equation*}
\|\vetx^{k+ \delta} - \vety \|^2 \leq \frac{1}{|\mathcal{M}|}\sum_{\ell \in \mathcal{M}} \| \vetx_{\ell}^{k+\delta} - \vety \|^2.
\end{equation*}

Using (\ref{missa1}) we obtain
\begin{align*} 
\| \vetx^{k+ \delta} - \vety \|^2 & \leq \frac{1}{|\mathcal{M}|} \sum_{\ell \in \mathcal{M}} \| \vetx^{k+ \delta-1} - \lambda_{k+ \delta -1}(\gks + \eks) - \vety \|^2 \\
  & = \frac{1}{|\mathcal{M}|} \sum_{\ell \in \mathcal{M}} [ \| \vetx^{k+ \delta-1} - \vety \|^2 - 2\lambda_{k+ \delta -1}(\gks+\eks)^{T}(\vetx^{k+ \delta-1}-\vety) \\
  & \quad + \lambda_{k+ \delta -1}^{2} \| \gks + \eks \|^2] \\
  &= \| \vetx^{k+ \delta-1} - \vety \|^2 - 2\lambda_{k+ \delta -1} \frac{1}{|\mathcal{M}|} \sum_{\ell \in \mathcal{M}} ( \gks + \eks )^{T}(\vetx^{k+ \delta-1} - \vety) \\
  & \quad +\lambda_{k+ \delta -1}^{2} \frac{1}{|\mathcal{M}|} \sum_{\ell \in \mathcal{M}} \| \gks+\eks \|^2.
\end{align*}
Repeating the same procedure for $\| \vetx^{k+ \delta-1} - \vety \|^2$, next for $\| \vetx^{k+ \delta-2} - \vety \|^2, \dots$, we obtain, after $\delta$ steps

\begin{align*} & \ds \| \vetx^{k+ \delta} - \vety \|^2 \leq \dxy - \frac{2}{|\mathcal{M}|}\sum_{q=1}^{\delta} \lambda_{k+\delta-q} \\
& \quad \times \sum_{\ell \in \mathcal{M}} \left( \mathbf{g}_{s_{\ell}(k+\delta-q+1)}^{k+\delta-q} + \boldsymbol \epsilon_{s_{\ell}(k+\delta-q+1)}^{k+\delta-q+1} \right)^{T}(\vetx^{k+\delta-q} - \vety) \\
							  & \quad + \frac{1}{|\mathcal{M}|} \sum_{q=1}^{\delta} \lambda_{k+\delta-q}^{2} \sum_{\ell \in \mathcal{M}} \| \mathbf{g}_{s_{\ell}(k+\delta-q+1)}^{k+\delta-q} + \boldsymbol \epsilon_{s_{\ell}(k+\delta-q+1)}^{k+\delta-q+1} \|^2.
\end{align*}
Using the definition of subdifferential sets and remembering that $j(q) = \delta - q + 1$, we have
\begin{align*} \| \vetx^{k+ \delta} - \vety \|^2 & \leq \dxy - \frac{2}{|\mathcal{M}|} \sum_{\ell \in \mathcal{M}} \left[ \sum_{q=1}^{\delta} \lambda_{k+\delta-q} (f_{s_\ell(k+j(q))}(\vetx^{k+\delta-q}) - f_{s_\ell(k+j(q))}(\vety)) \right.\\				& \quad + \left. \sum_{q=1}^{\delta} \lambda_{k+\delta-q} \left(\eksq \right)^{T}(\vetx^{k+\delta-q} - \vety) \right] + \frac{1}{|\mathcal{M}|} \sum_{\ell \in \mathcal{M}} \sum_{q=1}^{\delta} \lambda_{k+\delta-q}^{2} \| \gksq+\eksq \|^2,
\end{align*}
and, alternatively
\begin{align*}
\| \vetx^{k+ \delta} - \vety \|^2 & \leq \dxy - \frac{2}{|\mathcal{M}|} \sum_{\ell \in \mathcal{M}} \left[ \sum_{q=1}^{\delta} \lambda_{k+\delta-q} (f_{s_\ell(k+j(q))}(\vetx^{k+\delta-q}) - f_{s_\ell(k+j(q))}(\vetx^{\delta n(\ceil{\frac{k}{\delta}})})) \right.\\
							     &  \quad + \left. \sum_{q=1}^{\delta} \lambda_{k+\delta-q}(f_{s_\ell(k+j(q))}(\vetx^{\delta n(\ceil{\frac{k}{\delta}})}) - f_{s_\ell(k+j(q))}(\vety))
							     + \sum_{q=1}^{\delta} \lambda_{k+\delta-q} \left(\eksq \right)^{T}(\vetx^{k+\delta-q} - \vety) \right] \\
							     & \quad \quad + \frac{1}{|\mathcal{M}|} \sum_{\ell \in \mathcal{M}} \sum_{q=1}^{\delta} \lambda_{k+\delta-q}^{2} \| \gksq+\eksq \|^2.
\end{align*}
By Assumption \ref{assump1} (and using (\ref{subgrad_boundedness})), we have
\begin{align*} \| \vetx^{k+ \delta} - \vety \|^2 & \leq \dxy - \frac{2}{|\mathcal{M}|} \sum_{\ell \in \mathcal{M}} \left[ \sum_{q=1}^{\delta} \lambda_{k+\delta-q} (f_{s_\ell(k+j(q))}(\vetx^{k+\delta-q}) - f_{s_\ell(k+j(q))}(\vetx^{\delta n(\ceil{\frac{k}{\delta}})})) \right. \\
							     & \quad + \left. \sum_{q=1}^{\delta} \lambda_{k+\delta-q}(f_{s_\ell(k+j(q))}(\vetx^{\delta n(\ceil{\frac{k}{\delta}})}) - f_{s_\ell(k+j(q))}(\vety)) + \sum_{q=1}^{\delta} \lambda_{k+\delta-q}\left(\eksq \right)^{T}(\vetx^{k+\delta-q} - \vety) \right] \\
							     & \quad \quad + \frac{1}{|\mathcal{M}|} \sum_{\ell \in \mathcal{M}} \sum_{q=1}^{\delta} \lambda_{k+\delta-q}^{2} (C^{2} + 2C \| \eksq \| + \| \eksq \|^2).
\end{align*}
By analyzing the previous inequality over the iterates $\delta k$ we have
\begin{align*} \| \vetx^{\delta (k+1)} - \vety \|^2 & \leq \| \vetx^{\delta k} - \vety \|^2 - \frac{2}{|\mathcal{M}|} \sum_{\ell \in \mathcal{M}} \left[  \sum_{q=1}^{\delta} \lambda_{\delta (k+1)-q} (f_{s_\ell(\delta k + j(q))}(\vetx^{\delta (k+1)-q}) - f_{s_\ell(\delta k + j(q))}(\vetx^{\delta n(k)})) \right. \\
							     & \quad + \left.  \sum_{q=1}^{\delta} \lambda_{\delta (k+1)-q} (f_{s_\ell(\delta k + j(q))}(\vetx^{\delta n(k)}) - f_{s_\ell(\delta k + j(q))}(\vety)) \right. \\
							     & \quad \left. +  \sum_{q=1}^{\delta} \lambda_{\delta (k+1)-q} \left(\eksdq \right)^{T}(\vetx^{\delta (k+1)-q} - \vety) \right] \\
							     & \quad \quad + \frac{1}{|\mathcal{M}|} \sum_{\ell \in \mathcal{M}}  \sum_{q=1}^{\delta} \lambda_{\delta (k+1)-q}^{2} (C^{2} + 2C \| \eksdq \| + \| \eksdq \|^2).
\end{align*}
By equation (\ref{stepsize1}) and taking conditional expectations with respect to $\mathcal{G}_{\delta n(k)}$, we have
\begin{align} \label{aux1}
\EE[\| \vetx^{\delta (k+1)} - \vety \|^2 | \mathcal{G}_{\delta n(k)}] & \leq \EE[ \| \vetx^{\delta k} - \vety \|^2 | \mathcal{G}_{\delta n(k)}] -  \frac{2\lambda_{\delta k}}{|\mathcal{M}|} \sum_{\ell \in \mathcal{M}} \sum_{q=1}^{\delta} \EE[f_{s_{\ell}(\delta k + j(q))}(\vetx^{\delta (k+1)-q}) \\
									& \qquad - f_{s_{\ell}(\delta k + j(q))}(\vetx^{\delta n(k)}) | \mathcal{G}_{\delta n(k)} ] \nonumber \\
									& \quad-  \frac{2\lambda_{\delta k}}{|\mathcal{M}|} \sum_{\ell \in \mathcal{M}} \sum_{q=1}^{\delta} \EE[f_{s_{\ell}(\delta k + j(q))}(\vetx^{\delta n(k)}) - f_{s_{\ell}(\delta k + j(q))}(\vety) | \mathcal{G}_{\delta n(k)} ] \nonumber\\
									& \quad -  \frac{2\lambda_{\delta k}}{|\mathcal{M}|} \sum_{\ell \in \mathcal{M}} \sum_{q=1}^{\delta}  \EE\left[\left( \eksdq \right)^{T}(\vetx^{\delta (k+1)-q}-\vety) | \mathcal{G}_{\delta n(k)} \right] \nonumber\\
									&\quad +  \frac{\lambda_{\delta k}^{2}}{|\mathcal{M}|} \sum_{\ell \in \mathcal{M}} \sum_{q=1}^{\delta}( C_{\ell}^{2} + 2C_{\ell} \EE[\| \eksdq \| | \mathcal{G}_{\delta n(k)}] + \EE[\| \eksdq \|^2 | \mathcal{G}_{\delta n(k)}]) \nonumber.
\end{align}
For each $q=1, \dots, \delta$ define $\tilde{\vetg}_{s_{\ell}(\delta k + j(q))} \in \partial f_{s_{\ell}(\delta k + j(q))}(\vetx^{\delta n(k)})$. Then, using Cauchy-Schwarz inequality and Assumption
\ref{assump1}, we can estimate the second term from the right-hand side of (\ref{aux1}), for each $\ell \in \mathcal{M}$ and $q=1,\dots, \delta$, as
\begin{align} \label{aux2}  \EE[f_{s_{\ell}(\delta k + j(q))}(\vetx^{\delta (k+1)-q}) - f_{s_{\ell}(\delta k + j(q))}(\vetx^{\delta n(k)}) | \mathcal{G}_{\delta n(k)} ] & \geq \EE[-(\tilde{\vetg}_{s_{\ell}(\delta k + j(q))})^{T}(\vetx^{\delta n(k)} - \vetx^{\delta (k+1)-q}) | \mathcal{G}_{\delta n(k)}] \nonumber \\
												   &\geq \EE[-\| \tilde{\vetg}_{s_{\ell}(\delta k + j(q))} \| \| \vetx^{\delta n(k)} - \vetx^{\delta (k+1)-q} \| | \mathcal{G}_{\delta n(k)}] \nonumber \\
												   &\geq -C \EE[\| \vetx^{\delta n(k)} - \vetx^{\delta (k+1)-q}\| | \mathcal{G}_{\delta n(k)}].
\end{align}
Now, we need to find an estimate for $\EE[\| \vetx^{\delta n(k)} - \vetx^{\delta (k+1)-q} \| | \mathcal{G}_{\delta n(k)}]$. To this end, we use (\ref{missa2})
and non-expansivity of the projection to obtain, for each $\ell \in \mathcal{M}$ and $q=1,\dots, \delta$
\begin{align*} \EE[\|\vetx^{\delta n(k)} - \vetx^{\delta (k+1)-q} \| | \mathcal{G}_{\delta n(k)}] & \leq \EE\left[ \sum_{p=\delta n(k)}^{\delta (k+1)-q-1}\| \vetx^{p+1} - \vetx^p \| | \mathcal{G}_{\delta n(k)} \right] \\
									   &= \sum_{p=\delta n(k)}^{\delta (k+1)-q-1} \EE \left[  \left\| \PX \left( \frac{1}{|\mathcal{M}|}\sum_{\ell \in \mathcal{M}} \vetx_{\ell}^{p+1}\right) - \vetx^p \right\| | \mathcal{G}_{\delta n(k)} \right] \\
									   &\leq \sum_{p=\delta n(k)}^{\delta (k+1)-q-1} \EE \left[  \left\| \frac{1}{|\mathcal{M}|}\sum_{\ell \in \mathcal{M}} \vetx_{\ell}^{p+1} - \vetx^p \right\| | \mathcal{G}_{\delta n(k)} \right] \\
									   &\leq \sum_{p=\delta n(k)}^{\delta (k+1)-q-1} \EE \left[ \frac{1}{|\mathcal{M}|} \sum_{\ell \in \mathcal{M}} \| \vetx_{\ell}^{p+1} - \vetx^p \| | \mathcal{G}_{\delta n(k)} \right].
\end{align*}
The law of iterated expectations, equation (\ref{missa1}) and Assumptions \ref{assump1}-\ref{assump2} provide
\begin{align} \label{aux3} \EE [ \| \vetx^{\delta n(k)} - \vetx^{\delta (k+1)-q} \| | \mathcal{G}_{\delta n(k)}] & \leq \sum_{p=\delta n(k)}^{\delta (k+1)-q-1} \lambda_{p} \frac{1}{|\mathcal{M}|} \sum_{\ell \in \mathcal{M}} \EE [ \| \vetg_{s_{\ell}(p+1)}^{p} \| + \| \boldsymbol \epsilon_{s_{\ell}(p+1)}^{p+1} \| | \mathcal{G}_{\delta n(k)} ] \nonumber \\
									   &\leq \sum_{p=\delta n(k)}^{\delta (k+1)-q-1} \lambda_{p} \frac{1}{|\mathcal{M}|} \sum_{\ell \in \mathcal{M}} ( \EE [ \| \vetg_{s_{\ell}(p+1)}^{r} \| | \mathcal{G}_{\delta n(k)}] + \EE [ \EE [ \| \boldsymbol \epsilon_{s_{\ell}(p+1)}^{p+1} \| | G_p] | \mathcal{G}_{\delta n(k)}]) \nonumber \\
									   &\leq \sum_{p=\delta n(k)}^{\delta (k+1)-q-1} \lambda_{p} \frac{1}{|\mathcal{M}|} \sum_{\ell \in \mathcal{M}}(C + \nu_{p+1}),
\end{align}
which holds with probability $1$.
Combining (\ref{aux2}) and (\ref{aux3}), we can rewrite (\ref{aux1}) as follows,
\begin{align} \label{aux4} \EE[\| \vetx^{\delta (k+1)} - \vety \|^2 | \mathcal{G}_{\delta n(k)}] & \leq \EE[ \| \vetx^{\delta k} - \vety \|^2 | \mathcal{G}_{\delta n(k)}] +  \frac{2\lambda_{\delta k}}{|\mathcal{M}|} \sum_{\ell \in \mathcal{M}} \sum_{q=1}^{\delta} \sum_{p=\delta n(k)}^{\delta (k+1)-q-1} \lambda_{p}(C + \nu_{p+1}) \nonumber\\
									&\quad -  \frac{2\lambda_{\delta k}}{|\mathcal{M}|} \sum_{\ell \in \mathcal{M}} \sum_{q=1}^{\delta} \EE[f_{s_{\ell}(\delta k + j(q))}(\vetx^{\delta n(k)}) - f_{s_{\ell}(\delta k + j(q))}(\vety) | \mathcal{G}_{\delta n(k)} ] \nonumber\\
									& \quad -  \frac{2\lambda_{\delta k}}{|\mathcal{M}|} \sum_{\ell \in \mathcal{M}} \sum_{q=1}^{\delta}  \EE\left[\left( \eksdq \right)^{T}(\vetx^{\delta (k+1)-q}-\vety) | \mathcal{G}_{\delta n(k)} \right] \nonumber\\
									& \quad +  \frac{\lambda_{\delta k}^{2}}{|\mathcal{M}|} \sum_{\ell \in \mathcal{M}} \sum_{q=1}^{\delta}( \nu_{\delta k + j(q)} + C )^2, \qquad \mbox{a.s.}
\end{align}
where in the last term we use Assumption \ref{assump2}. Since $\mathcal{G}_{\delta n(k)} \subset \mathcal{G}_{\delta (k+1)-q} \subset \mathcal{G}_{\delta k + j(q)}$ for any $q=1, \dots, \delta$, we use again the law of iterated expectations and
Assumption \ref{assump2} to estimate for each $\ell \in \mathcal{M}$ and $q=1,\dots, \delta$
\begin{align*} - \EE \left[ \left( \eksq \right)^{T}(\vetx^{k+\delta-q}-\vety) | \mathcal{G}_{\delta n(k)} \right] & = - \EE [ \EE \left[ \left( \eksdq \right)^{T}(\vetx^{\delta (k+1)-q} - \vety) | \mathcal{G}_{\delta (k+1)-q} \right] | \mathcal{G}_{\delta n(k)}] \\
									&= \EE [ \EE \left[ \eksdq | \mathcal{G}_{\delta (k+1)-q} \right]^{T}(\vety - \vetx^{\delta (k+1)-q}) | \mathcal{G}_{\delta n(k)} ] \\
									&\leq \EE [ \| \EE [ \eksdq | \mathcal{G}_{\delta (k+1)-q} ] \| \| \vetx^{\delta (k+1)-q} - \vety \| | \mathcal{G}_{\delta n(k)} ] \\
									&\leq \nu_{\delta k + j(q)} \EE[ \| \vetx^{\delta (k+1)-q} - \vety \| | \mathcal{G}_{\delta n(k)} ] \qquad \mbox{a.s.}
\end{align*}
and replacing in (\ref{aux4}), we have with probability $1$
\begin{align*} \EE[\| \vetx^{\delta (k+1)} - \vety \|^2 | \mathcal{G}_{\delta n(k)}] & \leq \EE[ \| \vetx^{\delta k} - \vety \|^2 | \mathcal{G}_{\delta n(k)}] +  2\lambda_{\delta k} \sum_{q=1}^{\delta} \sum_{p=\delta n(k)}^{\delta (k+1)-q-1} \lambda_{p}(C + \nu_{p+1}) \\
									&\quad -  \frac{2\lambda_{\delta k}}{|\mathcal{M}|} \sum_{\ell \in \mathcal{M}} \sum_{q=1}^{\delta} \EE[f_{s_{\ell}(\delta k + j(q))}(\vetx^{\delta n(k)}) - f_{s_{\ell}(\delta k + j(q))}(\vety) | \mathcal{G}_{\delta n(k)} ] \\
									&\quad +  2\lambda_{\delta k} \sum_{q=1}^{\delta} \nu_{\delta k + j(q)} \EE[ \| \vetx^{\delta (k+1)-q} - \vety \| | \mathcal{G}_{\delta n(k)} ] +  \lambda_{\delta k}^{2} \sum_{q=1}^{\delta}( \nu_{\delta k + j(q)} + C )^2.
\end{align*}
Taking expectations, we have
\begin{align} \label{aux5} \EE[\| \vetx^{\delta (k+1)} - \vety \|^2 ]  &\leq \EE[ \| \vetx^{\delta k} - \vety \|^2 ] +  2\lambda_{\delta k} \sum_{q=1}^{\delta} \sum_{p=\delta n(k)}^{\delta (k+1)-q-1} \lambda_{p}(C + \nu_{p+1}) \nonumber\\
									&\quad -  \frac{2\lambda_{\delta k}}{|\mathcal{M}|} \sum_{\ell \in \mathcal{M}} \sum_{q=1}^{\delta} \EE[f_{s_{\ell}(\delta k + j(q))}(\vetx^{\delta n(k)}) - f_{s_{\ell}(\delta k + j(q))}(\vety) ] \nonumber\\
									&\quad +  2\lambda_{\delta k} \sum_{q=1}^{\delta} \nu_{\delta k + j(q)} \EE[ \| \vetx^{\delta (k+1)-q} - \vety \| ]  +  \lambda_{\delta k}^{2} \sum_{q=1}^{\delta}( \nu_{\delta k + j(q)} + C )^2.
\end{align}
As regards $ \EE[f_{s_{\ell}(\delta k + j(q))}(\vetx^{\delta n(k)}) - f_{s_{\ell}(\delta k + j(q))}(\vety)]$, we have for each $\ell \in \mathcal{M}$ and $q=1,\dots, \delta$,
\begin{align} \label{aux6} \EE[f_{s_{\ell}(\delta k + j(q))}(\vetx^{\delta n(k)}) - f_{s_{\ell}(\delta k + j(q))}(\vety) ] &= \sum_{v=1}^{N} \EE[ ( f_{s_{\ell}(\delta k + j(q))}(\vetx^{\delta n(k)}) - f_{s_{\ell}(\delta k + j(q))}(\vety) ) 1_{\{s_\ell(\delta n(k) )\in\mathcal R_v \}} ] \nonumber\\ 
& \quad + \EE[ ( f_{s_{\ell}(\delta k + j(q))}(\vetx^{\delta n(k)}) - f_{s_{\ell}(\delta k + j(q))}(\vety) ) 1_{\{s_\ell(\delta n(k) )\in\mathcal T \}} ] \nonumber \\
& = \sum_{v=1}^{N} \sum_{i \in \mathcal{R}_v} \EE[ ( f_{i}(\vetx^{\delta n(k)}) - f_{i}(\vety) ) 1_{\{s_\ell(\delta k + j(q))=i, \, s_\ell(\delta n(k) )\in\mathcal R_v \}} ] \nonumber \\ 
& \quad + \EE[ ( f_{s_{\ell}(\delta k + j(q))}(\vetx^{\delta n(k)}) - f_{s_{\ell}(\delta k + j(q))}(\vety) ) 1_{\{s_\ell(\delta n(k) )\in\mathcal T \}} ].
\end{align}
Let us denote
\begin{align*} o_{\ell,q}(k) &= \EE[ ( f_{s_{\ell}(\delta k + j(q))}(\vetx^{\delta n(k)}) - f_{s_{\ell}(\delta k + j(q))}(\vety) ) 1_{\{s_\ell(\delta n(k) )\in\mathcal T \}} ]
\end{align*}
and rewrite \eqref{aux6} as
\begin{align*} \EE[f_{s_{\ell}(\delta k + j(q))}(\vetx^{\delta n(k)}) - f_{s_{\ell}(\delta k + j(q))}(\vety) ]
& = \sum_{v=1}^{N} \sum_{i,h \in \mathcal{R}_v} \EE[ ( f_i(\vetx^{\delta n(k)}) - f_i(\vety) ) 
\\& \qquad \times 1_{\{s_\ell(\delta k + j(q))=i, \, s_\ell(\delta n(k) )\in\mathcal R_v, \,
s_\ell(\delta n(k)+1)=h \}} ]+o_{\ell,q}(k).
\end{align*}
Let us also denote
\begin{align*}
\tau_{\ell, \, q}(k, \, i, \, h) &= \text{Prob}( s_\ell(\delta k + j(q)) = i , s_\ell(\delta n(k) +1) = h,s_\ell(\delta n(k) )\in\mathcal R_v) \quad \mbox{and}  \\
\tilde{\tau}_{\ell, \, q}(k, \, i, \, h) &= \text{Prob}( s_\ell(\delta k + j(q)) = i \,     |  \, s_\ell(\delta n(k) + 1) = h, \, s_\ell(\delta n(k) )\in\mathcal R_v),
\end{align*}
and then we have for each $\ell \in \mathcal{M}$ and $q=1,\dots, \delta$,
\begin{align*} \EE[f_{s_{\ell}(\delta k + j(q))}(\vetx^{\delta n(k)}) - f_{s_{\ell}(\delta k + j(q))}(\vety) ]
- o_{\ell,q}(k) & = \sum_{v=1}^{N} \sum_{i,h \in \mathcal{R}_v} \EE[ ( f_i(\vetx^{\delta n(k)}) - f_i(\vety) ) \\
& \qquad \times 1_{\{s_\ell(\delta k + j(q))=i, \, s_\ell(\delta n(k) + 1)= h, s_\ell(\delta n(k) )\in\mathcal R_v \}}  ] \\
& = \sum_{v=1}^{N} \sum_{i,h \in \mathcal{R}_v}  \EE[ ( f_i(\vetx^{\delta n(k)}) - f_i(\vety) ) \, | \, s_\ell(\delta k + j(q)) = i, \\
& \qquad s_\ell(\delta n(k) + 1) = h, s_\ell(\delta n(k) )\in\mathcal R_v ] \tau_{\ell, \, q}(k, \, i, \, h)  \\
& = \sum_{v=1}^{N} \sum_{i,h \in \mathcal{R}_v}  \EE[ ( f_i(\vetx^{\delta n(k)}) - f_i(\vety) ) \, | \, s_\ell(\delta n(k)+1) = h, \\
& \qquad s_\ell(\delta n(k) )\in\mathcal R_v ] \tau_{\ell, \, q}(k, \, i, \, h),
\end{align*}
where the last equality follows from $\delta k + j(q) \geq \delta n(k) + 1$ for all $q=1,\dots, \delta$.
Since $\tau_{\ell, \, q}(k, \, i, \, h) = \tilde{\tau}_{\ell, \, q}(k, \, i, \, h) \text{Prob}(s_\ell(\delta n(k) + 1) = h, s_\ell(\delta n(k) )\in\mathcal R_v)$, then we can write the preceding equation as
\begin{align*}  \EE[f_{s_{\ell}(\delta k + j(q))}(\vetx^{\delta n(k)}) - f_{s_{\ell}(\delta k + j(q))}(\vety) ]
- o_{\ell,q}(k) & = \sum_{v=1}^{N} \sum_{i,h \in \mathcal{R}_v}  \EE[ ( f_i(\vetx^{\delta n(k)}) - f_i(\vety) ) \\
& \qquad \times 1_{ \{ s_\ell( \delta n(k) + 1) = h, s_\ell(\delta n(k) )\in\mathcal R_v \}} ] \tilde{\tau}_{\ell, \, q}(k, \, i, \, h) \\
& = \sum_{v=1}^{N} \sum_{i,h \in \mathcal{R}_v}  \EE[ ( f_i(\vetx^{\delta n(k)}) - f_i(\vety) ) \\
& \qquad \times 1_{\{s_\ell(\delta n(k)+1) = h, s_\ell(\delta n(k) )\in\mathcal R_v \}} ] [ P^{\delta (k-n(k)) + j(q) - 1} ]_{h i}.
\end{align*}
Now we can use Lemma \ref{lem-eval-n1-2} to ensure that there are positive scalars $\eta$ and $\beta$ such that, for each $\ell \in \mathcal{M}$
\begin{align*} &\sum_{q=1}^{\delta} \Big( \EE [ f_{s_\ell(\delta k + j(q))}(\vetx^{\delta n(k)}) - f_{s_\ell(\delta k + j(q))}(\vety)] - o_{\ell,q}(k) \Big) \\
& = \sum_{v=1}^{N} \sum_{q=1}^{\delta} \sum_{i,h \in \mathcal{R}_v}  \EE[ ( f_i(\vetx^{\delta n(k)}) - f_i(\vety) ) 1_{\{s_\ell(\delta n(k)+1) = h, s_\ell(\delta n(k) )\in\mathcal R_v\}} ] [ P^{\delta (k-n(k)) + j(q) - 1} ]_{h i}  \\
& \geq
\sum_{v=1}^{N} \Big\{
- \eta e^{-\beta (k-n(k))} + \delta  \sum_{i \in \mathcal R_v} \EE \Big[ ( f_i(\vetx^{\delta n(k)}) - f_i(\vety) )  1_{ s_\ell(\delta n(k) )\in\mathcal R_v } \Big] [\boldsymbol \pi_{\delta,v}]_i\Big\} \\
& \geq - N \eta e^{-\beta (k-n(k))} + \delta  \sum_{v=1}^{N} \sum_{i \in \mathcal R_v} \EE \Big[ ( f_i(\vetx^{\delta n(k)}) - f_i(\vety) ) 1_{\{ s_\ell(\delta n(k) )\in\mathcal R_v \}} \Big] [ \boldsymbol \pi_{\delta,v}]_i \\ 
& \geq - N \eta e^{-\beta (k-n(k))} + \delta \sum_{v=1}^{N} \sum_{i \in \mathcal{R}_v} \Big( \text{Prob}(s_\ell(\delta n(k))\in\mathcal{R}_v) [ \boldsymbol \pi_{\delta,v}]_i \EE[f_{i}(\vetx^{\delta n(k)})- f_{i}(\vety) | s_\ell(\delta n(k))\in\mathcal{R}_v] \Big).
\end{align*}
Taking expectations we have
\begin{align*} & \EE \Bigg[  \sum_{q=1}^{\delta} \Big( \EE [ f_{s_\ell(\delta k + j(q))}(\vetx^{\delta n(k)}) - f_{s_\ell(\delta k + j(q))}(\vety)] -  o_{\ell,q}(k) \Big) \Bigg] \\
\\ & \geq \EE \Bigg[ - N \eta e^{-\beta (k-n(k))} + \delta \sum_{v=1}^{N} \sum_{i \in \mathcal{R}_v} \Big( 
\text{Prob}(s_\ell(\delta n(k))\in\mathcal{R}_v) [ \boldsymbol \pi_{\delta,v}]_i
\EE[f_{i}(\vetx^{\delta n(k)})- f_{i}(\vety) | s_\ell(\delta n(k))\in\mathcal{R}_v] \Big) \Bigg] \\
& =  - N \eta e^{-\beta (k-n(k))} + \delta \sum_{v=1}^{N} \sum_{i \in \mathcal{R}_v} \Big(
\text{Prob}(s_\ell(\delta n(k))\in\mathcal{R}_v) [ \boldsymbol \pi_{\delta,v}]_i
\EE \Big[ \EE[f_{i}(\vetx^{\delta n(k)})- f_{i}(\vety) | s_\ell(\delta n(k))\in\mathcal{R}_v] \Big] \Big),
\end{align*}
leading to
\begin{align*} & \sum_{q=1}^{\delta} \Big( \EE [ f_{s_\ell(\delta k + j(q))}(\vetx^{\delta n(k)}) - f_{s_\ell(\delta k + j(q))}(\vety)] -  o_{\ell,q}(k) \Big) \\
& \geq  - N \eta e^{-\beta (k-n(k))} + \delta \sum_{v=1}^{N} \sum_{i \in \mathcal{R}_v} \Big(
\text{Prob}(s_\ell(\delta n(k))\in\mathcal{R}_v) [ \boldsymbol \pi_{\delta,v}]_i
\EE[f_{i}(\vetx^{\delta n(k)})- f_{i}(\vety)] \Big).
\end{align*}
For simplicity, let us write
\begin{align}\label{aux7} \frac{1}{|\mathcal{M}|} \sum_{\ell \in \mathcal{M}} \sum_{q=1}^{\delta}
\EE [ f_{s_\ell(\delta k + j(q))}(\vetx^{\delta n(k)}) - f_{s_\ell(\delta k + j(q))}(\vety)] 
& \geq - N \eta e^{-\beta (k-n(k))} + o(k) \\
& \quad + \delta \sum_{v=1}^{N} \sum_{i \in \mathcal{R}_v} \Bigg( \frac{1}{|\mathcal{M}|} \sum_{\ell \in \mathcal{M}} [\boldsymbol \pi_{\ell}^{\infty} ]_i  \nonumber
\EE[f_{i}(\vetx^{\delta n(k)})- f_{i}(\vety)]   \Bigg),  \nonumber
\end{align}
where
\begin{align} \label{aux8} o(k) & = \frac{1}{|\mathcal{M}|} \sum_{\ell \in \mathcal{M}} \left( \sum_{q=1}^{\delta} o_{\ell,q}(k) + \delta \sum_{v=1}^{N} \sum_{i \in \mathcal{R}_v} \Big[ (\text{Prob}(s_{\ell}(\delta n(k)) \in \mathcal{R}_v)[ \boldsymbol \pi_{\delta,v}]_i - [\boldsymbol \pi_{\ell}^{\infty} ]_i) \EE[f_{i}(\vetx^{\delta n(k)})- f_{i}(\vety)] \Big] \right).
\end{align}
By using (\ref{weights}) and (\ref{main-problem})-\textbf{(i)}, we can rewrite \eqref{aux7} as
\begin{align} \label{eq-aux-eduardo01} \frac{1}{|\mathcal{M}|} \sum_{\ell \in \mathcal{M}} \sum_{q=1}^{\delta}
\EE [ f_{s_\ell(\delta k + j(q))}(\vetx^{\delta n(k)}) - f_{s_\ell(\delta k + j(q))}(\vety)] & \geq - N \eta e^{-\beta (k-n(k))}  + o(k) \nonumber \\
& \quad + \delta \sum_{v=1}^{N} \sum_{i \in \mathcal{R}_v} \left( \frac{1}{|\mathcal{M}|} \sum_{\ell \in \mathcal{M}} [\boldsymbol \pi_{\ell}^{\infty} ]_i \right) \EE[f_{i}(\vetx^{\delta n(k)})- f_{i}(\vety) ]
\nonumber \\
& = - N \eta e^{-\beta (k-n(k))}  + o(k) \nonumber \\
& \quad + \delta \sum_{v=1}^{N}  \sum_{i \in \mathcal{R}_v} \omega_i
\EE[f_{i}(\vetx^{\delta n(k)})- f_{i}(\vety) ] \nonumber \\
& = - N \eta e^{-\beta (k-n(k))}  + o(k) \nonumber \\
& \quad + \delta \sum_{v=1}^{N} \EE[ \sum_{i \in \mathcal{R}_v} \omega_i  (f_{i}(\vetx^{\delta n(k)})- f_{i}(\vety)) ] \nonumber \\
& = - N \eta e^{-\beta (k-n(k))} + \delta \EE[ f(\vetx^{\delta n(k)})- f(\vety) ] + o(k).
\end{align}
At this juncture we move back to the evaluation of the term $o(k)$.
The compactness of $\X$ and convexity of $f_i$ guarantee
the existence of $\bar M > 0$ such that
$|f_{s_{\ell}(\delta k + j(q))}(\vetx^{\delta n(k)}) - f_{s_{\ell}(\delta k + j(q))}(\vety)|\leq \bar M$ a.s.; this and statement (i) of Corollary \ref{cor-transient} yield
\begin{equation*}\begin{aligned} & \Big| \EE \Big[ (
f_{s_{\ell}(\delta k + j(q))}(\vetx^{\delta n(k)}) - f_{s_{\ell}(\delta k + j(q))}(\vety)
) 1_{\eliasc{\{}{}\{s_\ell(\delta n(k) )\in\mathcal T\}}  \Big] \Big|
\\& \leq
\EE \Big[\Big|  (
f_{s_{\ell}(\delta k + j(q))}(\vetx^{\delta n(k)}) - f_{s_{\ell}(\delta k + j(q))}(\vety)
) 1_{\eliasc{\{}{}\{s_\ell(\delta n(k) )\in\mathcal T\}}  \Big|\Big]
\\& \leq  \bar M \EE \Big[ 1_{\eliasc{\{}{}\{s_\ell(\delta n(k) )\in\mathcal T\}} \Big]
\leq \bar M \bar{\alpha} e^{- \bar{\beta} \eliasc{k}{n( k )}}
\end{aligned}
\end{equation*}
or, equivalently, $|o_{\ell,q}(k)|\leq \bar M\eliasc{}{\bar{\alpha}e^{- \bar{\beta} n( k )}}$,
so that $o_{\ell,q}(k) \geq -\bar M\eliasc{}{\bar{\alpha} e^{- \bar{\beta} n( k )}}$.
Next we substitute this inequality in \eqref{aux8}
and we use statement (ii) of Corollary \ref{cor-transient}
and the fact that the same $\bar M$ as above provides $|f_{i}(\vetx^{\delta n(k)})- f_{i}(\vety)|\leq \bar M$ a.s.,
and we obtain
\begin{align*} o(k) & \geq  \frac{1}{|\mathcal{M}|} \sum_{\ell \in \mathcal{M}} \left( \sum_{q=1}^{\delta} (-\bar M)\eliasc{}{\bar{\alpha} e^{- \bar{\beta} n( k )}} - \bar M \delta \sum_{v=1}^{N} \sum_{i \in \mathcal{R}_v} \Big[ \text{Prob}(s_{\ell}(\delta n(k)) \in \mathcal{R}_v)[ \boldsymbol \pi_{\delta,v}]_i - [\boldsymbol \pi_{\ell}^{\infty} ]_i \Big]\right) \\ 
& \geq \frac{-\delta \bar M}{|\mathcal{M}|} \sum_{\ell \in \mathcal{M}} 
\eliasc{\sum_{v=1}^{N} \sum_{i \in \mathcal{R}_v} \bar{\alpha} e^{- \bar{\beta} k}}%
{\left( \bar{\alpha} e^{- \bar{\beta} n( k )} + \sum_{v=1}^{N} \sum_{i \in \mathcal{R}_v} \bar{\alpha} e^{- \bar{\beta} n( k )} \right)} \\
& \geq  \frac{-\delta \bar M \eliasc{N |\mathcal{M}|}{|\mathcal{M}|(1 + m)}}{|\mathcal{M}|}  \bar{\alpha} e^{- \bar{\beta} \eliasc{k}{n( k )}}
\\
& =  -\delta \bar M (1 + m) \bar{\alpha} e^{- \bar{\beta} n( k )}
\end{align*}
and from \eqref{eq-aux-eduardo01},
\begin{align*} \frac{1}{|\mathcal{M}|} \sum_{\ell \in \mathcal{M}} \sum_{q=1}^{\delta}
\EE [ f_{s_\ell(\delta k + j(q))}(\vetx^{\delta n(k)}) - f_{s_\ell(\delta k + j(q))}(\vety)] & \geq - \varrho (e^{-\theta (k-n(k))} +  e^{- \theta n( k )}) + \delta \EE[ f(\vetx^{\delta n(k)})- f(\vety) ],
\end{align*}
where we defined
$\varrho =  \max(N \eta, \delta \bar M \eliasc{N}{( 1 + m )} \bar{\alpha} )$\eliasc{}{ and $\theta = \min( \beta, \bar{\beta} )$}.
The above and equation (\ref{aux5}) lead to \eqref{eq-lemma-main-estimate}, completing the proof.
\end{proof}

We can now approach the proof of our main result.
\begin{proof}[Proof of Theorem~\ref{teorema1}] Since $\X$ is compact and $f$ is convex over $\mathbb{R}^n$ (therefore, also continuous), the optimal set $\X^{\ast}$ is
nonempty, compact and convex. Define $n(0)=0$ and $n(k) = k + 1 - \ceil{k^\gamma}$ ($k \geq 1$) with $1-\xi < \gamma < 2 \xi - 1$.
Notice that $n(k) \leq k$ for all $k \geq 0$. Therefore, we can use Lemma \ref{lemma-main-estimate} with
this particular choice for $n(k)$. 
Taking $\vety \eliasc{= \PXast(\vetx^{\delta k})}{{}\in \X^*}$\eliasc{and using the relation
$d(\vetx^{\delta (k+1)}, \, \X^{\ast}) \leq \| \vetx^{\delta (k+1)} - \PXast(\vetx^{\delta k}) \|$}{}, Lemma \ref{lemma-main-estimate} provides
\begin{align}  \label{aux13} \EE[\| \vetx^{\delta (k+1)} - \vety \|^2 ] & \leq \EE[\| \vetx^{\delta k} - \vety \|^2 ] +  \frac{2a}{(k+1)^\xi} \sum_{q=1}^{\delta} \sum_{p=\delta n(k)}^{\delta (k+1)-q-1} \lambda_{p} \Big(C + \max_{p=\delta n(k), \dots, \delta (k+1) - q - 1}{\nu_{p+1}} \Big) \nonumber \\
									& \quad -  \frac{2 a \delta}{(k+1)^\xi}  \EE[ f(\vetx^{\delta n(k)}) - f^{\ast}] + \frac{2a \varrho \delta}{(k+1)^\xi} \Big(e^{- \eliasc{\bar{\beta}}{\theta} n(k)} + e^{- \eliasc{\beta}{\theta} (k - n(k))} \Big) \nonumber \\
									& \quad +  \frac{2a}{(k+1)^\xi} \sum_{q=1}^{\delta} \nu_{\delta k + j(q)} \EE[ \| \vetx^{\delta (k+1)-q} - \eliasc{\PXast(\vetx^{\delta k})}{\vety} \| ] +  \frac{a^2}{(k+1)^{2\xi}} \sum_{q=1}^{\delta}( \nu_{\delta k + j(q)} + C )^2.
\end{align}
Define $d_{\X} = \max_{\vetx, \, \vetz \in \X}{\| \vetx - \vetz \|}$ and
\begin{align*}
\rho_k &= \frac{2a(C + \nu)}{(k+1)^\xi} \sum_{q=1}^{\delta} \sum_{p=\delta n(k)}^{\delta (k+1)-q-1} \lambda_{p} + \frac{2a \varrho \delta}{(k+1)^\xi} \Big(e^{- \eliasc{\beta}{\theta} n(k)} + e^{- \eliasc{\beta}{\theta} (k - n(k))} \Big) +  \frac{2 a d_{\X}}{(k+1)^\xi}  \sum_{q=1}^{\delta} \nu_{\delta k + j(q)} \\ & \quad +  \frac{\delta a^2 ( \nu + C )^2}{(k+1)^{2\xi}}.
\end{align*}
We can rewrite the inequality (\ref{aux13}) as
\begin{align} \label{aux14} \EE[\| \vetx^{\delta (k+1)} - \vety \|^2 ] & \leq \EE[\| \vetx^{\delta k} - \vety \|^2 ]
 - \frac{2 a\delta}{(k+1)^\xi} \EE[ f(\vetx^{\delta n(k)}) - f^{\ast} ] + \rho_k.
\end{align}
We next show that $\sum_{k=0}^{\infty} \rho_k < \infty$. Notice that,
\begin{align*} \frac{1}{(k+1)^\xi}  \sum_{q=1}^{\delta} \sum_{p=\delta n(k)}^{\delta (k+1)-q-1} \lambda_{p} & = \frac{1}{(k+1)^\xi} \Big[  ( \overbrace{\lambda_{\delta n(k)} + \dots + \lambda_{\delta k + \delta - 2}}^{q=1} + \dots + \overbrace{\lambda_{\delta n(k)} + \dots + \lambda_{\delta k - 1}}^{q= \delta} ) \Big] \\
      & \leq \frac{1}{(k+1)^\xi} \Big[  ( \lambda_{\delta n(k) - 1} + \dots + \lambda_{\delta n(k) - 1} + \dots + \lambda_{\delta n(k) - 1} + \dots + \lambda_{\delta n(k) - 1} ) \Big],
\end{align*}
because $\delta \geq 1$ and $n(k) \leq k$. Due to equation (\ref{stepsize-theorem}), we have
\begin{equation*}
\lambda_{\delta n(k)-1} = \lambda_{\delta (n(k) - 1)} = \frac{a}{n(k)^\xi} = \frac{a}{(k+1 - \ceil{k^\gamma})^\xi}.
\end{equation*}
Therefore,
\begin{align} \label{aux15} \frac{1}{(k+1)^\xi} \sum_{q=1}^{\delta} \sum_{p=\delta n(k)}^{\delta (k+1)-q-1} \lambda_{p} & \leq \frac{1}{(k+1)^\xi} \frac{a}{(k+1 - \ceil{k^\gamma})^\xi} \sum_{q=1}^{\delta} (\delta(k-n(k)+1)-q-1) \nonumber \\
& < a \delta^2 \frac{\ceil{k^\gamma}}{(k+1)^\xi (k+1 - \ceil{k^\gamma})^\xi},
\end{align}
because $-q-1 < 0$. Notice that $\ceil{k^\gamma} - 1 < k^\gamma$ for all $k \geq 0$ and thus
\begin{align*}
\frac{\ceil{k^\gamma}}{(k+1)^\xi (k+1 - \ceil{k^\gamma})^\xi} & < \frac{k^\gamma + 1}{k^\xi (k - k^\gamma)^\xi} \\
&= \frac{k^\gamma + 1}{k^{2\xi} (1 - k^{\gamma-1})^\xi} \\
&\leq \frac{k^\gamma + 1}{k^{2\xi} (1 - 2^{\gamma-1})^\xi},
\end{align*}
holds for all $k \geq 2$ since $\gamma - 1 < 0$. Hence, for all $k \geq 2$ we have
\begin{align*}
\frac{\ceil{k^\gamma}}{(k+1)^\xi (k+1 - \ceil{k^\gamma})^\xi} & < \frac{1}{(1 - 2^{\gamma-1})^\xi} \frac{k^\gamma + 1}{k^{2\xi}} \\
&= \frac{1}{(1 - 2^{\gamma-1})^\xi} \left(\frac{1}{k^{2\xi - \gamma}} + \frac{1}{k^{2 \xi}}\right) \end{align*}
and together with inequality (\ref{aux15}) we conclude that
\begin{align*} \sum_{k=2}^{\infty} \left( \frac{1}{(k+1)^\xi} \sum_{q=1}^{\delta} \sum_{p=\delta n(k)}^{\delta (k+1)-q-1} \lambda_{p} \right) & < \frac{a \delta^2}{(1 - 2^{\gamma-1})^\xi} \left( \sum_{k=2}^{\infty} \frac{1}{k^{2\xi - \gamma}} + \sum_{k=2}^{\infty} \frac{1}{k^{2 \xi}}\right) < \infty,
\end{align*}
because $2\xi - \gamma > 1$. Furthermore, since $0 < \gamma < 1$, there is $\bar k > 0$ such that $k \geq \bar k \Rightarrow k+1 - \ceil{k^\gamma} \geq \ceil{k^\gamma}$ and thus
\begin{align*} 2a \varrho \delta \sum_{k=\bar{k}}^{\infty} \frac{1}{(k+1)^\xi} \Big(e^{-\eliasc{\beta}{\theta} n(k)} +  e^{- \eliasc{\beta}{\theta} (k - n(k)) )} \Big) & < 2a \varrho \delta \sum_{k=\bar{k}}^{\infty} \left( \frac{1}{k^\xi e^{\eliasc{\beta}{\theta} (k+1-\ceil{k^\gamma})} } + \frac{e^\eliasc{\beta}{\theta}}{k^\xi e^{\eliasc{\beta}{\theta} (\ceil{k^\gamma})}} \right) < \infty,
\end{align*}
because $e^x > x \,\, \forall x \in \mathbb{R}$, $\ceil{k^\gamma} > k^\gamma$ and $\xi + \gamma > 1$.
By assumption (\ref{assump_stepsize-moment}) and noticing that $2 \xi > 1$, we conclude that $\sum_{k=0}^{\infty} \rho_k < \infty$.
Denoting $a_k := \eliasc{d(\vetx^{\delta k}, \X^{\ast})}{\| \vetx^{\delta k} - \vety \|}^2$, $b_k \equiv 0$ and
$c_k := \frac{2a \delta}{(k+1)^\xi} (f(\vetx^{\delta n(k)}) - f^{\ast})$, inequality (\ref{aux14}) provides
\begin{align*}
\EE[a_{k+1}] \leq (1+ b_k) \EE[a_k] - \EE[c_k] + \EE[\rho_k].
\end{align*}
We can use Lemma \ref{lema2} (in the surely sense) to conclude that $\EE[a_k]$ converges to a non-negative scalar and
\begin{align*}
\sum_{k=0}^{\infty} \EE[c_k] = 2 a \delta \sum_{k=0}^{\infty} \frac{1}{(k+1)^\xi} \EE[f(\vetx^{\delta n(k)}) - f^{\ast}] < \infty.
\end{align*}
Since $\sum_{k=0}^{\infty} 1 / (k+1)^\xi = \infty$ and $f(\vetx^{\delta n(k)}) \geq f^{\ast}$, we have
\begin{align}\label{eq-comenario-eduardo}
\lim_{k \to \infty} \inf \EE[f(\vetx^{\delta n(k)})] = f^{\ast}.
\end{align}
The function $f(\vetx)$ is bounded on $\X$ (because $f$ is convex over $\R^n$ and $\X$ is bounded), thus
\begin{align*} \lim_{k \to \infty} \inf f(\vetx^{\delta n(k)}) \geq f^{\ast},
\end{align*}
and from Fatou's lemma we obtain
\begin{align*} \EE[ \lim_{k \to \infty} \inf f(\vetx^{\delta n(k)})] \leq \lim_{k \to \infty} \inf \EE[f(\vetx^{\delta n(k)})] = f^{\ast}.
\end{align*}
The two preceding inequalities imply that $\lim_{k \to \infty} \inf f(\vetx^{\delta n(k)}) = f^{\ast}$ with probability $1$. This relation, together with
the continuity of $f$ and boundedness of $\X$, implies that $\lim_{k \to \infty} \inf d(\vetx^{\delta n(k)}, \X^{\ast}) = 0$ with probability $1$.
\end{proof}
\begin{proof}[Proof of Theorem~\ref{teorema2}]
Notice that Lemma~\ref{lemma-main-estimate} with $\lambda_k \equiv \lambda > 0$ gives
\begin{align*} \frac1\lambda\EE[\| \vetx^{\delta (k+1)} - \vety \|^2]  & \leq \frac1\lambda\EE[ \| \vetx^{\delta k} - \vety \|^2] +  2 \lambda\sum_{q=1}^{\delta} \sum_{p=\delta n(k)}^{\delta (k+1)-q-1} (C + \nu_{p+1}) -  2 \delta \EE[ f(\vetx^{\delta n(k)}) - f(\vety)] \nonumber \\
& \quad +  2 \sum_{q=1}^{\delta} \nu_{\delta k + j(q)} \EE[ \| \vetx^{\delta (k+1)-q} - \vety \| ] + \lambda \sum_{q=1}^{\delta}( \nu_{\delta k + j(q)} + C )^2 \nonumber \\
& \qquad + 2 \varrho \delta (e^{- \theta n(k)} + e^{-\theta (k-n(k))}).
\end{align*}
Recall that $\nu \ge \nu_k$ for all $k \in \mathbb N$ and that, from Assumption~\ref{assump1}, we know that there is $M \ge \| \vetx - \vety \|$ for all $\vetx$, $\vety \in \X$. Then, letting $n( k ) = k - \zeta$ (we will soon specify $\zeta$), we have
\begin{align*} \frac1{2\lambda\delta}\EE[\| \vetx^{\delta (k+1)} - \vety \|^2]  & \leq \frac1{2\lambda\delta}\EE[ \| \vetx^{\delta k} - \vety \|^2] +  \lambda\zeta\delta(C + \nu) + \varrho e^{-\theta \zeta} - \EE[ f(\vetx^{\delta n(k)}) - f(\vety)] + \varrho e^{- \theta n(k)} \\ & \quad +  M\nu \nonumber + \frac\lambda2 ( \nu + C )^2.
\end{align*}

We select $\zeta$ as a function of the scalars $\lambda$, $\delta$, $C$, $\nu$, $\varrho$, and $\theta$ in a way to (almost) minimize the error estimate as follows. Let us define
$$
\phi( \zeta ) := \lambda\zeta \delta( \nu + C ) + \varrho e^{-\theta \zeta}.
$$
If $\varrho > 0$, then $\phi$ is strictly convex and has the minimizer
$$
\zeta^* = -\frac{\ln \frac{\lambda \delta( \nu + C )}{\theta\varrho}}\theta.
$$
Notice that if $\lambda$ is small enough, then $\zeta^* > 0$. For $\zeta^* > 0$, we may instead use $\zeta = \lceil \zeta^* \rceil$. Notice that $\zeta^* \le \zeta < \zeta^* + 1$, therefore, $\phi( \zeta ) < \phi( \zeta^* + 1 )$. Thus,
\begin{align*} \frac1{2\lambda\delta}\EE[\| \vetx^{\delta (k+1)} - \vety \|^2]  & \leq \frac1{2\lambda\delta}\EE[ \| \vetx^{\delta k} - \vety \|^2] +  \lambda(\zeta^* + 1)\delta(C + \nu) + \varrho e^{-\theta (\zeta^* + 1)} - \EE[ f(\vetx^{\delta n(k)}) - f(\vety)] \nonumber \\
& \quad +  \varrho e^{- \theta n(k)} + M\nu \nonumber + \frac\lambda2 ( \nu + C )^2,
\end{align*}
that is
\begin{equation}\label{eq:elias_aux_estimate}
\begin{split} \frac1{2\lambda\delta}\EE[\| \vetx^{\delta (k+1)} - \vety \|^2] & \leq \frac1{2\lambda\delta}\EE[ \| \vetx^{\delta k} - \vety \|^2] + \varrho e^{- \theta n(k)} + \frac{\delta( \nu + C )}\theta\lambda\ln\frac1\lambda + A\lambda \\ 
& \quad + M\nu - \EE[ f(\vetx^{\delta n(k)}) - f(\vety)],
\end{split}
\end{equation}
where $A = \frac{\delta( \nu + C )}\theta\left(- \ln\frac{\delta( \nu + C )}{\theta\varrho} + \theta + e^{-\theta}\right) + \frac{( \nu + C )^2}2$.

Now, assume for contradiction that
\begin{equation*}
      \liminf_{k \to \infty}\EE[ f(\vetx^{\delta k}) - f^{\ast}] > \frac{\delta( \nu + C )}\theta\lambda\ln\frac1\lambda + A\lambda + M\nu.
\end{equation*}
Then, for some $\epsilon > 0$ and $k_0 \in \mathbb N$, we have that for all $k \ge k_0$
\begin{equation}\label{eq:elias-large-deviation}
\EE[ f(\vetx^{\delta n( k )}) - f^{\ast}] > \frac{\delta( \nu + C )}\theta\lambda\ln\frac1\lambda + A\lambda + M\nu + \epsilon.
\end{equation}
Moreover, we can also assume that $k_0$ is large enough such that $\varrho e^{-\theta n( k )} \leq \epsilon / 2$ for $k \ge k_0$. Putting this and \eqref{eq:elias-large-deviation} in~\eqref{eq:elias_aux_estimate}, we get
\begin{equation}
\frac1{2\lambda\delta}\EE[\| \vetx^{\delta (k+1)} - \vety \|^2] \leq \frac1{2\lambda\delta}\EE[ \| \vetx^{\delta k} - \vety \|^2] - \frac\epsilon2.
\end{equation}
Iterating and rearranging we get
\begin{equation}
\EE[\| \vetx^{\delta ( k + n )} - \vety \|^2] \leq \EE[ \| \vetx^{\delta k} - \vety \|^2] - n\lambda\delta\epsilon,
\end{equation}
which is a contradiction because the right-hand side is negative for large enough $n$.
If $\lambda$ is not sufficiently small and we have $\zeta^* \le 0$, the result is proven using $\zeta = 1$ and following a similar reasoning.
\end{proof}%
\eduardo{
\begin{remark} \label{rem-comparing-with-sundhar2009incremental}
If we consider ergodic Markov chains with uniform limiting distribution
(i.e., $\omega_i = 1/m$ for all $i \in \II$) and $|\mathcal{M}| = 1$, then Theorem~\ref{teorema1} recovers Theorem 4.3 of \cite{sundhar2009incremental}.
However, the evaluations behind our main result are different from the
ones in \cite{sundhar2009incremental}.
The use of more general Markov chains makes several passages of the convergence analysis considerably more complex in comparison with the analysis of the method with ergodic Markov chains. We had to consider some extra terms
in the proof of Lemma \ref{lemma-main-estimate} (when compared e.g. with \cite[Lemma 4.2]{sundhar2009incremental} in the
simplified scenario of $\omega_i = 1/m$ and $|\mathcal{M}| = 1$), see equation~\eqref{aux6}, that concerns the inclusion of recurrent periodic and transient states and all subsequent equations. Corollary \ref{cor-transient} and Lemma \ref{lema4} given in the Appendix and all related evaluations are not required in simplified scenarios.
\end{remark}}

\section{Experimental results} \label{sec.4}
In this section, we report the numerical results of a simple example to illustrate a possible situation we can handle when using Algorithm \ref{missa}: we consider minimizing the $\ell_1$-norm of the residual associated to the linear system $A \vetx = \mathbf{b}$. We choose a matrix $A = (a_{i,j})$ ($7 \times 20$) with a high degree of sparsity and a feasible set $\X$ given by box constraints. The nonzero entries of $A$ are given by:
{\small\begin{align*}
a_{1,2} &= 0.5; \, a_{1,3} = 0.1; \, a_{1,4} = 0.2; a_{1,14} = 0.25; \, a_{1,15} = 0.1; \\
a_{2,6} &= 0.4; \, a_{2,7} = 0.15; \, a_{2,12} = 0.3; \, a_{2,16} = 0.45;  \\
a_{2,19} &= 0.1; \, a_{2,20} = 0.2; \\
a_{3,13} &= 0.02; \, a_{3,14} = 0.06; \\
a_{4,1} &= 0.12; \, a_{4,2} = 0.21; \, a_{4,3} = 0.3; a_{4,7} = 0.5; \, a_{4,13} = 0.4; \\
\,\,a_{4,14} &= 0.1; \, a_{4,15} = 0.18; \, a_{4,19} = 0.1; \, a_{4,20} = 0.14; \\
a_{5,1} &= 0.8; \, a_{5,2} = 0.4; \, a_{5,8} = 1.2; a_{5,9} = 1.0; \, a_{5,10} = 0.85; \\
\,\, a_{5,17} &= 0.4; \, a_{5,18} = 0.7; \, a_{5,19} = 0.1; \\
a_{6,2} &= 0.25; \, a_{6,3} = 0.34; \, a_{6,8} = 0.45; a_{6,9} = 0.35; \\
a_{6,13} &= 0.18; a_{6,14} = 0.22; \\
a_{7,13} &= 0.05; \, a_{7,14} = 0.08.
\end{align*}}
The feasible set is such that
$$\vetx \in \X \Leftrightarrow l_j \leq x_j \leq u_j, \,\, \forall j = 1, \dots 20,$$
with $l_j$, $x_j$ and $u_j$ denoting the $j$th component of $\mathbf{l}$, $\vetx$ and $\mathbf{u}$, respectively. The vectors $\mathbf{l}$ and $\mathbf{u}$ were chosen as
\begin{align*}
\mathbf{l} &= [-1, \, -0.5, \, -1.5, \, -1.3, \, 0, \,\, 0.1, \,\, 0.3, \, -0.2, \, -1.0, \,\, 0, \\
& -0.25, \, -0.1, \,\, 0.3, \,\, 0.1, \,\, 0, \, -1.1, \,\, 0.35, \,\, 0.15, \, 0, \,\, -0.45]^{T}
\end{align*}
and
\begin{align*}
\mathbf{u} &= [2.0, \, 1.5, \, 2.3, \, 3.0, \, 2.0, \,\, 1.8, \,\, 2.25, \, 1.7, \, 1.5, \,\, 2.0, \\
& 2.8, \, 1.75, \,\, 2.35, \,\, 1.95, \,\, 2.0, \, 1.0, \,\, 2.5, \,\, 1.35, \, 2.0, \,\, 3.0]^{T}.
\end{align*}

In order for the optimal set $\X^{\ast}$ to be nonempty, we set a \rafaelb{vector $\vety \in \X$ and compute $\mathbf{b} = A \vety$.}
We have constructed the transition probability matrix
$$ \small P=\begin{bmatrix}
0 & 0 & 0.2 & 0.8 & 0 & 0 & 0 \\
0 & 0 & 0.15 & 0.85 & 0 & 0 & 0 \\
0.4 & 0.6 & 0 & 0 & 0 & 0 & 0 \\
0.5 & 0.5 & 0 & 0 & 0 & 0 & 0 \\
0 & 0 & 0 & 0 & 0 & 0.8 & 0.2 \\
0 & 0 & 0 & 0 & 0.8 & 0 & 0.2 \\
0 & 0 & 0 & 0 & 0.6 & 0.4 & 0
\end{bmatrix}.
$$
Following the notation of (\ref{trans_prob_matrix}), we can note that $P = \text{diag}(P_1, \, P_2)$ and we have two \eduardo{irreducible} sets of recurrent states: $\mathcal{R}_1 = \{1,2,3,4\}$ and $\mathcal{R}_2 = \{ 5,6,7\}$. States in $\mathcal{R}_1$ are of period $\delta_1 = 2$, while states in $\mathcal{R}_2$ are aperiodic ($\delta_2 = 1$) and thus, $P^2$ is the transition probability matrix of an aperiodic Markov chain.
We test Algorithm \ref{missa} with two Markov chains, $s_1$ and $s_2$, with initial distributions
\begin{align*}
\boldsymbol \pi_{1}^{0} &= [1, \, 0, \, 0, \, 0, \, 0, \, 0, \, 0]^{T} \quad \text{and} \\
\boldsymbol \pi_{2}^{0} &= [0, \, 0, \, 0, \, 0, \, 1, \, 0, \, 0]^{T},
\end{align*}
respectively. By using equation (\ref{weights}) we compute the weights $\omega_i$:
\begin{align*}
\boldsymbol \omega & := [\omega_1, \dots, \omega_7]^{T} \\
& \approx [0.121, \, 0.129, \, 0.043, \,  0.206, \, 0.213, \, 0.203, \, 0.083]^{T},
\end{align*}
in such a way we have a correspondence between the entries of $\boldsymbol \omega$
and the rows of $A$: the larger is the norm of the $i$-th row of $A$, the
larger is $\omega_i$.

Define $W = \text{diag}(\omega_1, \dots, \omega_7)$.
We have $A \vetx^{\ast} = \mathbf{b}$ iff $WA \vetx^{\ast} = W\mathbf{b}$. Thus, we can consider the equivalent optimization problem:
\begin{equation*} \begin{array}{c} \ds \vetx \in \arg \min f(\vetx) = \| W A \vetx - W \mathbf{b} \|_1
\\
\mbox{s.t.} \quad \vetx \in \X, \end{array} \end{equation*}
which has the same form that problem (\ref{main-problem}). In fact,
$$f(\vetx) = \sum_{i=1}^{7} \omega_i |\mathbf{a}_{i}^{T} \vetx - b_i|,$$
where $\mathbf{a}_{i}^{T}$ represents the $i$th row of $A$. The subgradients $\vetg_i$ of $f_i(\vetx) = |\mathbf{a}_{i}^{T} \vetx - b_i|$ can be computed by the rule:
\begin{equation*}
 \vetg_i = \left\{ \begin{array}{ccc}
                    \mathbf{a}_{i}, & \text{if} & \mathbf{a}_{i}^{T} \vetx - b_i > 0, \\
                    - \mathbf{a}_{i}, & \text{if} & \mathbf{a}_{i}^{T} \vetx - b_i < 0, \\
                    \mathbf{0}, && \text{otherwise},
                   \end{array} \right.
\end{equation*}
because $A^{T} \mathbf{sgn}(A \vetx - \mathbf{b}) \in \partial \| A \vetx - \mathbf{b}\|_1$ (see, e.g., \cite{hiriart93I}), where $\mathbf{sgn}$ is the sign function.

We compare the performance of Algorithm \ref{missa} with three other methods: the incremental (cyclic) subgradient method, Markov randomized incremental subgradient method and incremental (randomized) subgradient method (\cite{nedic01},  \cite{sundhar2009incremental} and \cite{johansson2010}).
\eduardo{All of them are} particular cases of the Algorithm \ref{missa} when $|\mathcal{M}| = 1$ (i.e., when we use just one Markov chain), with specific transition probability matrices. Note that these methods must consider subgradients of $h_i(\vetx) := \omega_i f_i(\vetx)$. For ease of notation, we label the methods -- Algorithm \ref{missa} (MISSA), Markov randomized incremental subgradient method, incremental (cyclic) subgradient method and incremental (randomized) subgradient method -- by $M1$, $M2$, $M3$ and $M4$ respectively. \rafaelc{The transition probability matrix for $M3$, denoted by $P_{M3}$, is given by $[P_{M3}]_{i,i+1}=1$, $i=1,\dots,6$, $[P_{M3}]_{7,1}=1$ and $[P_{M3}]_{i,j}=0$ for all other entries. The transition probability matrix for $M4$, denoted by $P_{M4}$, is given by $[P_{M4}]_{i,j} = 1/7$ for all $i,j \in \{1,\dots,7\}$.}
\rafaelc{Moreover, the initial distributions for $M3$ and $M4$ are $[1, \, 0, \dots, 0]^{T}$ and $[1/7, \, 1/7, \dots, 1/7]^{T}$, respectively.} For the method $M2$, the transition probability matrix (or all matrices $P(k)$ in \cite{sundhar2009incremental})\rafaelc{, denoted by $P_{M2}$,} must satisfy some assumptions that we do not consider \eduardo{in} Algorithm \ref{missa}, among them, irreducibility, aperiodicity and uniform limiting distribution. In our tests, we consider one of the suggested rules in \cite{sundhar2009incremental}
(equal probability scheme):
$$[P_{M2}]_{i,j} = \left\{ \begin{array}{ccc}
                           1/7, & \text{if} & j \neq i \, \text{and} \, j \in N_i, \\
                           1 - \frac{|N_i|}{7}, & \text{if} & j=i, \\
                           0, & & \text{otherwise},
                           \end{array} \right.
$$
where $N_i \subset \{1, \dots, 7 \}$ is the set of neighbors of an agent $i$. We set $N_1 = \{2,3\}$, $N_2 = \{1,3,7\}$, $N_3 = \{1,2,6\}$, $N_4 = \{5,6\}$, $N_5 = \{4\}$, $N_6 = \{3,4,7\}$ and $N_7 = \{2,6\}$, besides an initial distribution \eduardo{equals} to $[0, \, 0, \, 0, \, 0, \, 1, \, 0, \, 0]^{T}$, aiming to get higher probabilities of reaching agent $5$ (the row of $A$ with larger norm) in the initial iterations.

Regarding the stochastic errors $\boldsymbol \epsilon_{i}^{k}$ of the subgradients for the agent $i$, we tested 
six  
different possibilities in order to study the effect of the Assumption \ref{assump2} and (\ref{assump_stepsize-moment})  (see similar assumptions for $M2$ and $M3$ in \cite{sundhar2009incremental}) on the performance of the methods. 
\rafaelc{In what follows, $\mathcal U(0,1)$ denotes the uniform distribution in the interval $[0, 1]$ and $\mathcal N(0,1)$ the normal distribution with zero mean and unit variance. In all cases, $[\boldsymbol \epsilon_i^k]_j$ and $[\boldsymbol \epsilon_{i'}^{k'}]_j$, $i\not=i'$, $k\not=k'$, $j=1,\dots,20$, are independent random variables with distributions given in Table~\ref{table1}.}
\begin{table}[htbp!]
\centering
\setlength{\arrayrulewidth}{2\arrayrulewidth}
\setlength{\belowcaptionskip}{10pt}
\begin{tabular}{|c|c|}
\hline
 Tests  & Description (for any $k > 0$ and $j=1, \dots, 20$) \\
\hline
Test $1$ & $\boldsymbol \epsilon_{i}^{k} = [0,0,\dots,0]^T$ for all $i \in \II$. \\
Test $2$ &  $[\boldsymbol \epsilon_{i}^{k}]_j \sim \mathcal U(0,k^{-1})$ for all $i \in \II$. \\
Test $3$ & $[\boldsymbol \epsilon_{i}^{k}]_j \sim 0.1 \mathcal{U}(0,1)$ for all $i \in \II$.  \\
Test $4$ &  $[\boldsymbol \epsilon_{i}^{k}]_j \sim 0.01 \mathcal{U}(0,1)$ for all $i \in \II$. \\
Test $5$ &  $[\boldsymbol \epsilon_{i}^{k}]_j \sim 0.1 \mathcal{N}(0,1)$ for all $i \in \II$. \\
Test $6$ &  $[\boldsymbol \epsilon_{i}^{k}]_j \sim 0.01 \mathcal{N}(0,1)$ for all $i \in \II$. \\
\hline
\end{tabular}
\caption{Description of the stochastic errors used in the tests.} \label{table1}
\end{table}

\rafaelb{We initially present results using a diminishing stepsize rule. The initial guess was chosen as $\PX([0,0,\dots,0]^T)$ for all methods. We run all methods during $10{,}000{,}000$ iterations.} The stepsizes for the algorithms $M2$, $M3$ and $M4$ was chosen as in Theorem 4.3 in \cite{sundhar2009incremental} and are similar to the sequence (\ref{stepsize-theorem}) adopted for $M1$ (with the exception of the periodicity $\delta$). The parameters $a$ and $\xi$ were tuned for each method, with $a \in \{1.0, 1.5, 2.0, 2.5, 3.0\}$ and $\xi \in \{0.667, 0.7, 0.8, 0.9, 1\}$ \rafaelb{and running the methods for test $1$ (without random errors).} Table~\ref{stepsizes_param} shows the parameters used in all tests for each method.
\begin{table}[htbp!]
\centering
\setlength{\arrayrulewidth}{2\arrayrulewidth}
\setlength{\belowcaptionskip}{15pt}
\begin{tabular}{|c|c|c|}
\hline
 Method  & $a$ & $\xi$ \\
\hline
$M1$ & $2.0$ & $0.7$ \\
$M2$ &  $2.0$ & $0.7$ \\
$M3$ & $2.5$ & $0.667$ \\
$M4$ &  $2.5$ & $0.667$ \\
\hline
\end{tabular}
\caption{Parameters used in the stepsizes of the methods.} \label{stepsizes_param}
\end{table}

All tests were run on a computer with an Intel(R) Core(TM) i3-4005U CPU @ 1.70GHz $\times$ 4 processor and $4$GB of RAM. We run $M1$ in parallel, following the procedure of allocating a thread for each Markov chain.
It is natural to expect the CPU time per iteration to be greater for $M1$ than for $M2$, $M3$ and $M4$. Although the subiterations in~\eqref{missa1} are computed in parallel for $M1$ (in our implementation, we use the \textsc{std::thread} library (C$++$11) and the code was compiled by the g$++$ compiler of GCC -- version 7.5.0) we need to compute, as we can see in~\eqref{missa2}, the average of the subiterations before calculating the projection. The average is not required for the other methods, as they only use one Markov chain, i.e., $|\mathcal{M}| = 1$. For example, in test $2$, the CPU time per iteration for $M1$ is about $5.2$e$-05$ seconds, while for $M2$, $M3$ and $M4$ it is approximately $3.4$e$-05$ seconds.

\begin{figure*}[ht!]
\begin{center}
\begin{tabular}{ccc}
\includegraphics[scale=0.37]{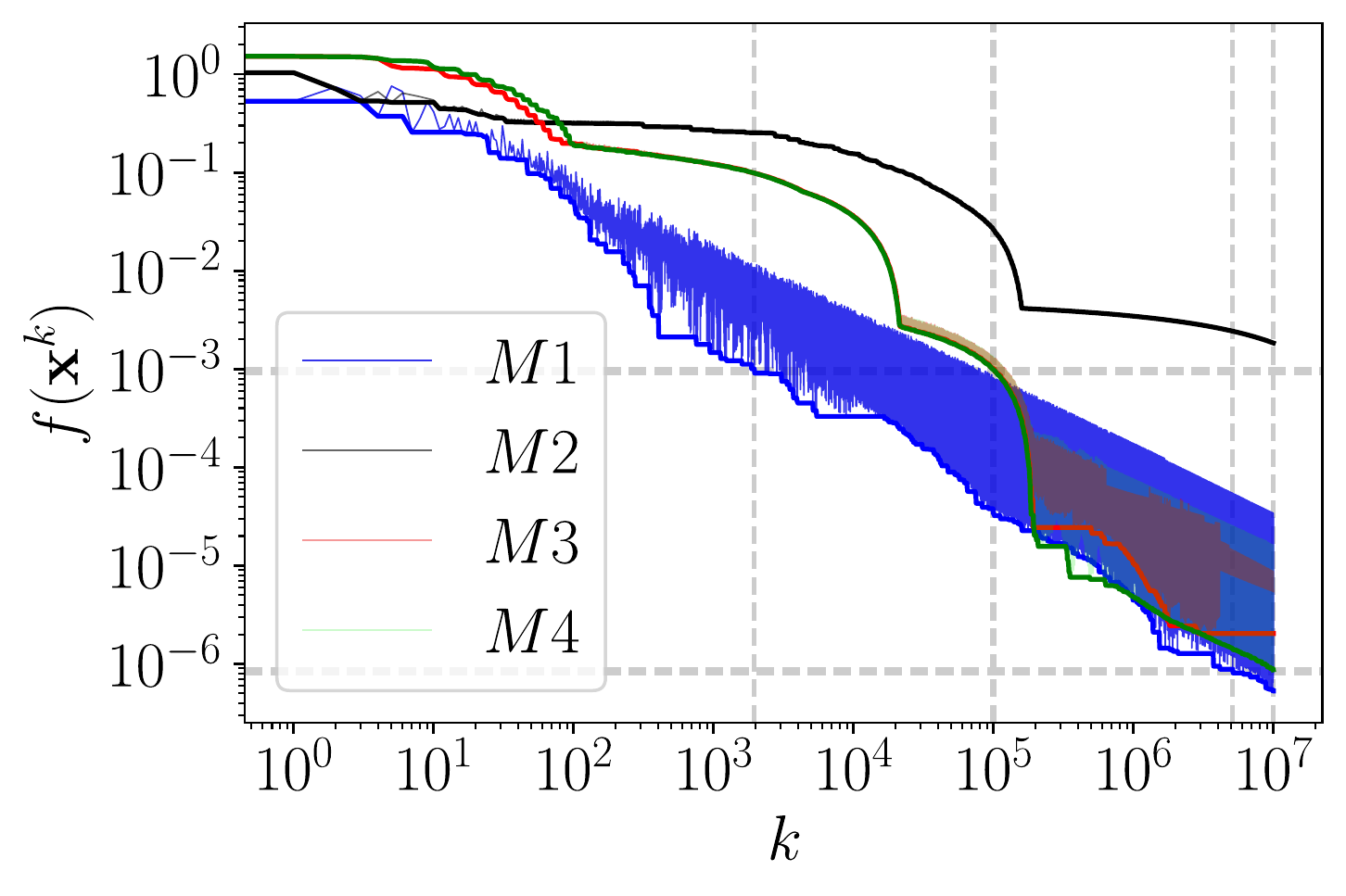} &
\includegraphics[scale=0.37]{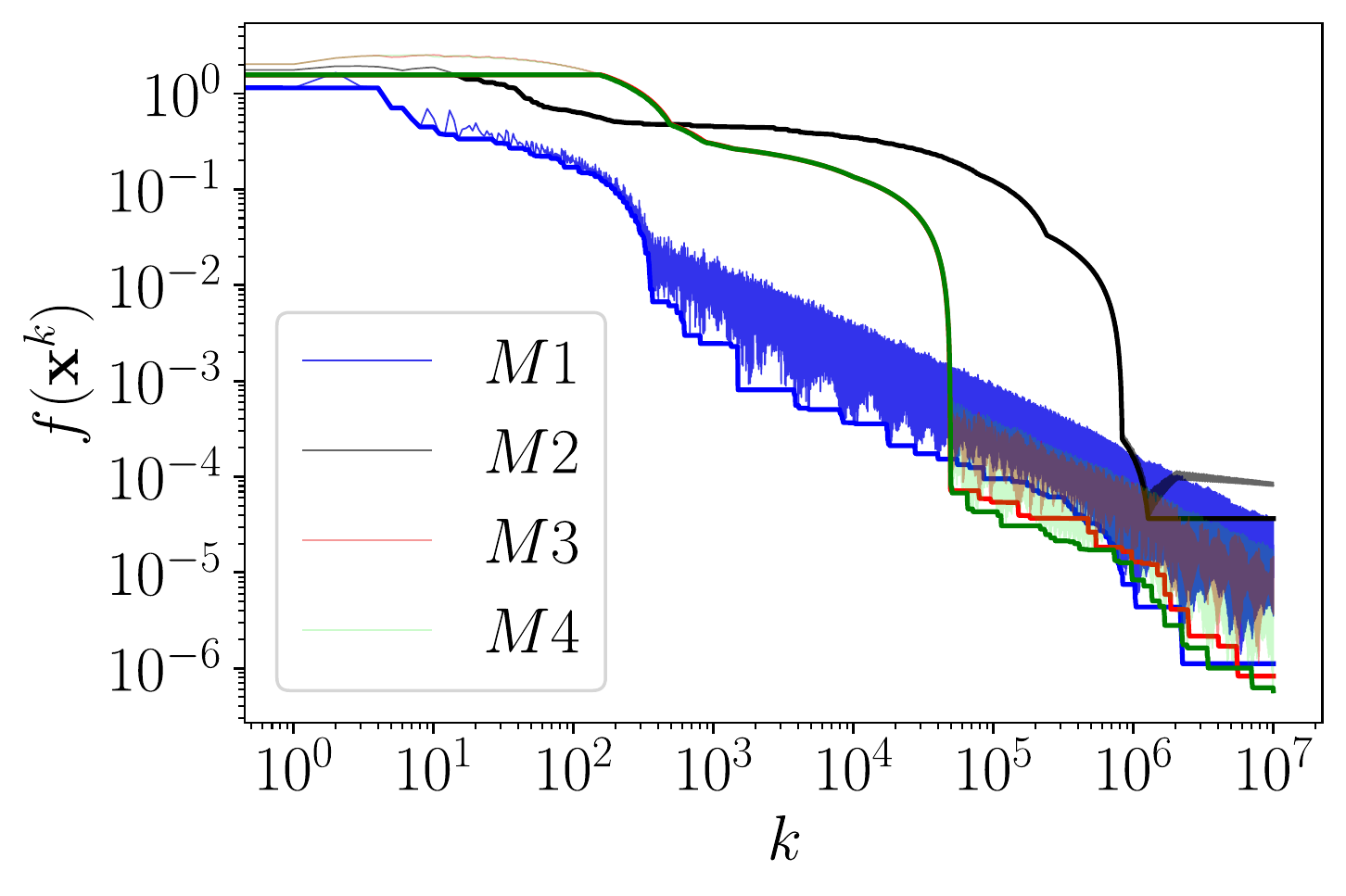} &
\includegraphics[scale=0.37]{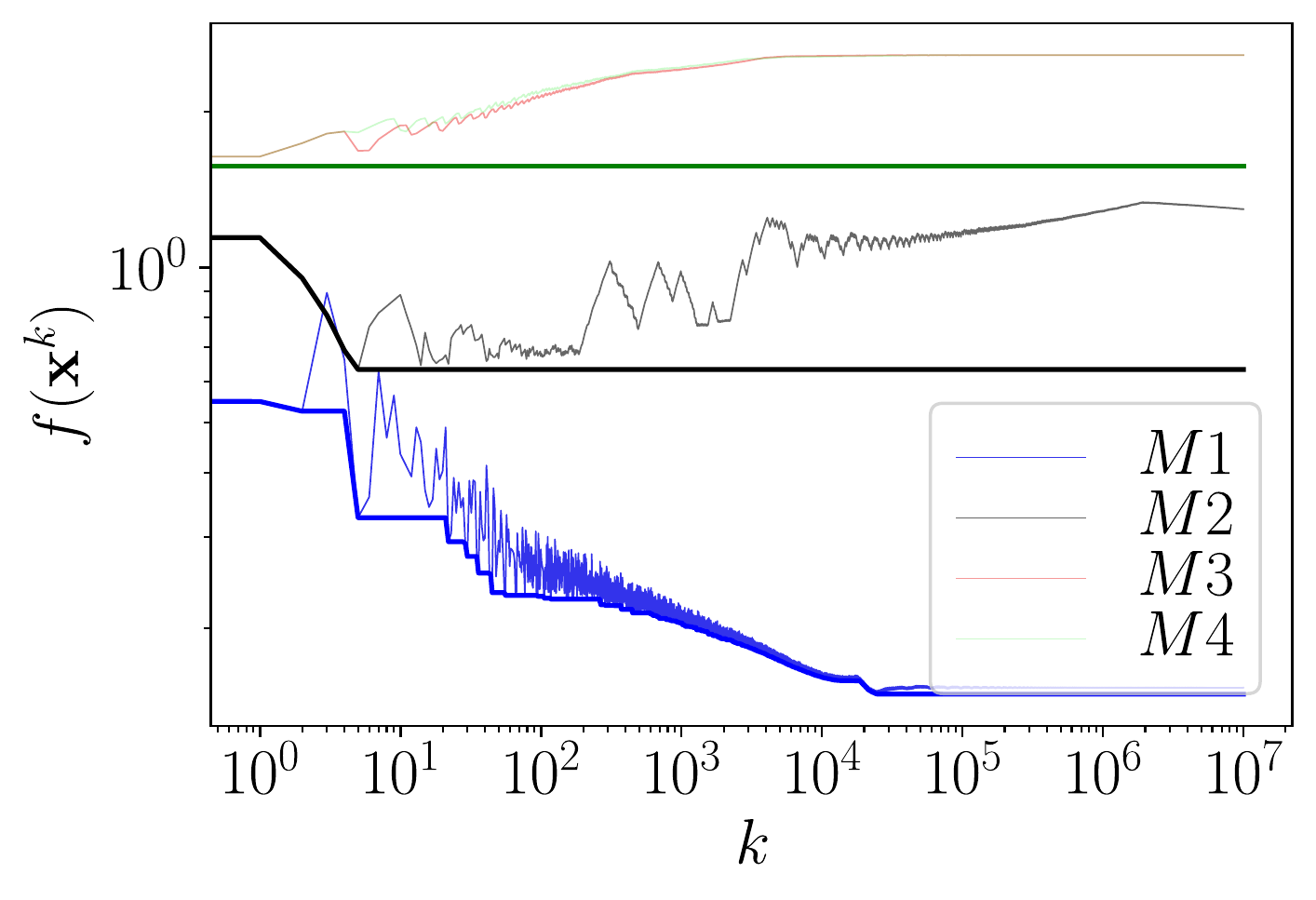} \\
(a) Test $1$ & (b) Test $2$ & (c) Test $3$ \\
\includegraphics[scale=0.37]{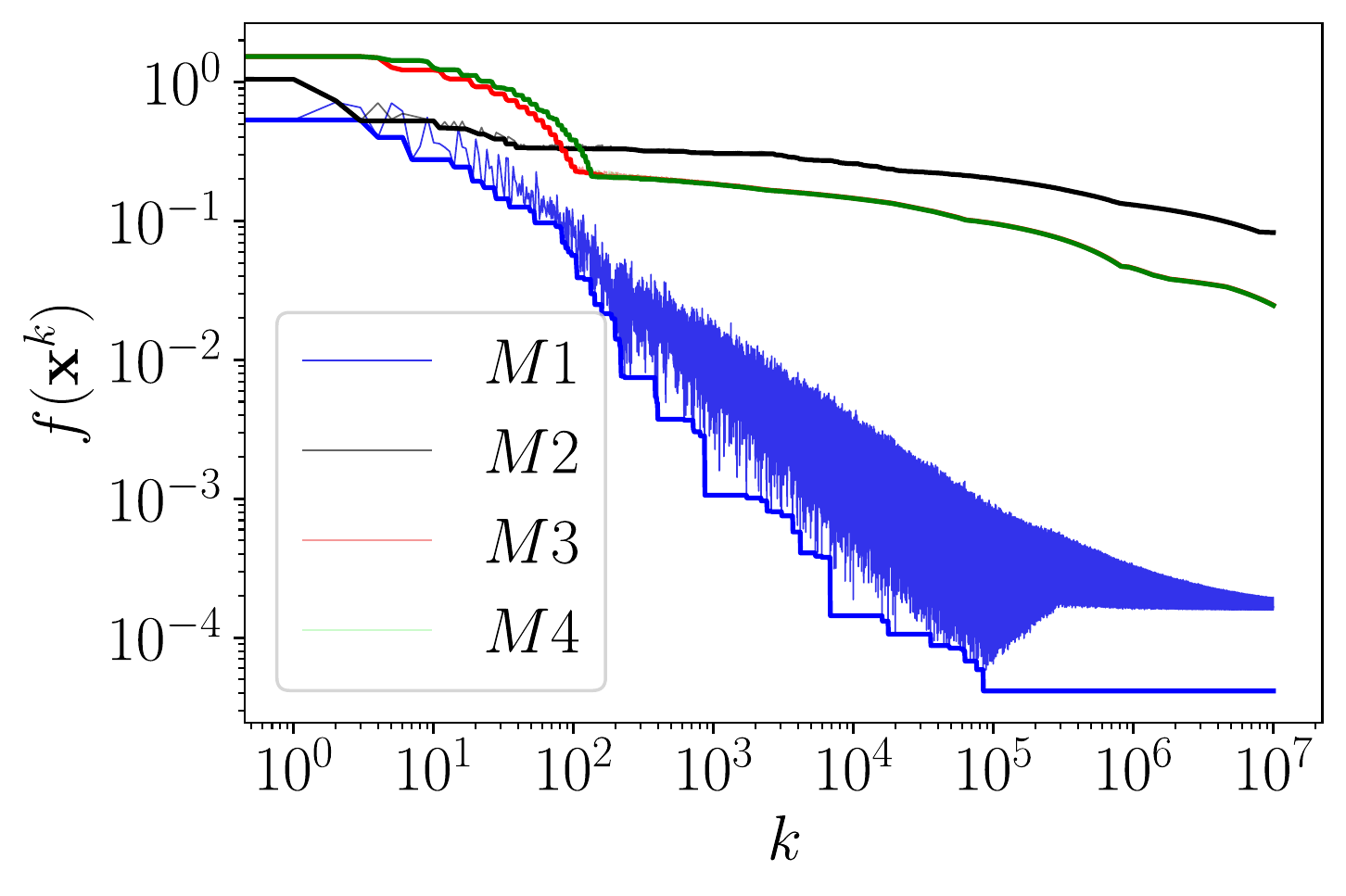} &
\includegraphics[scale=0.37]{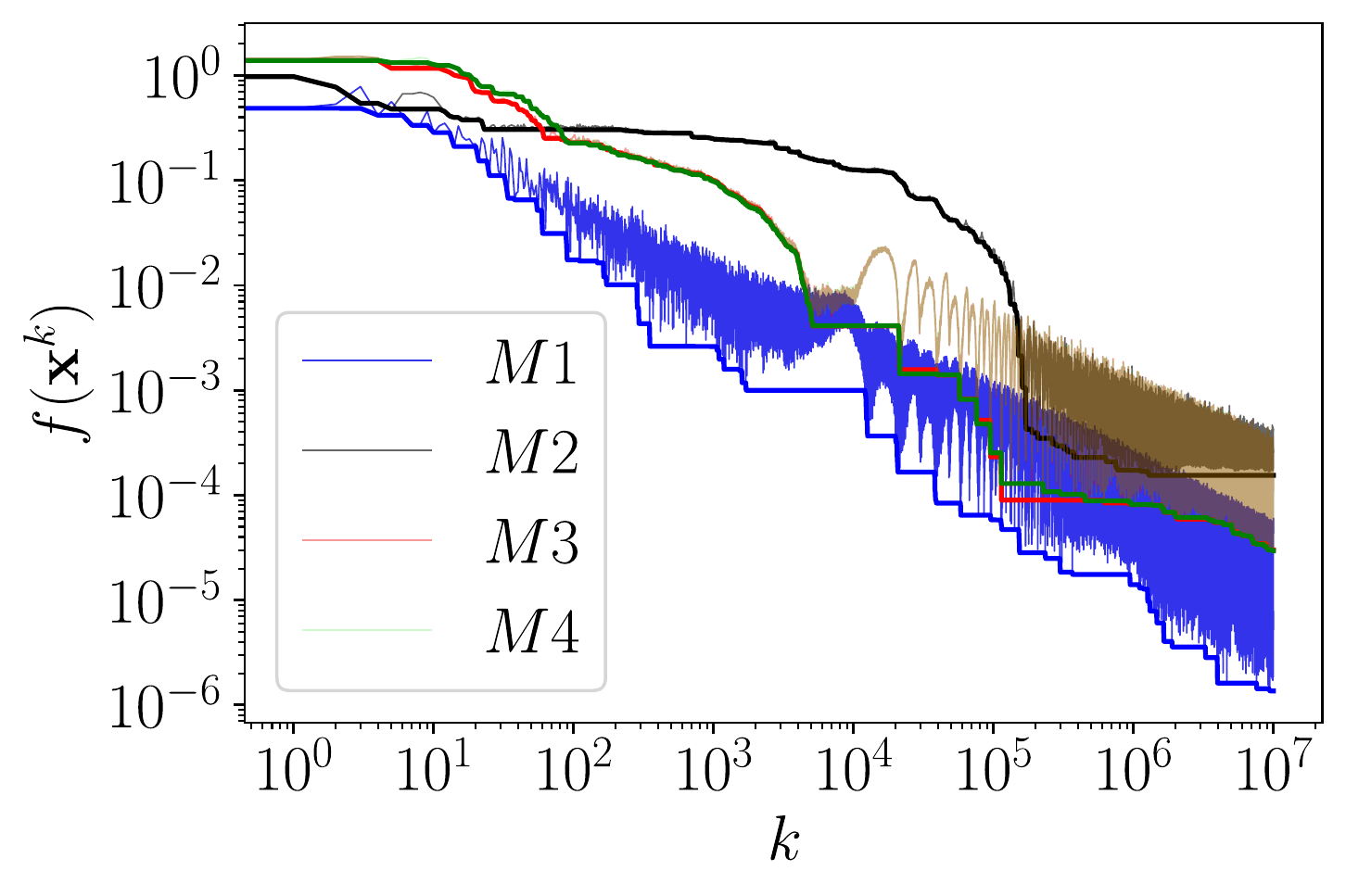} &
\includegraphics[scale=0.37]{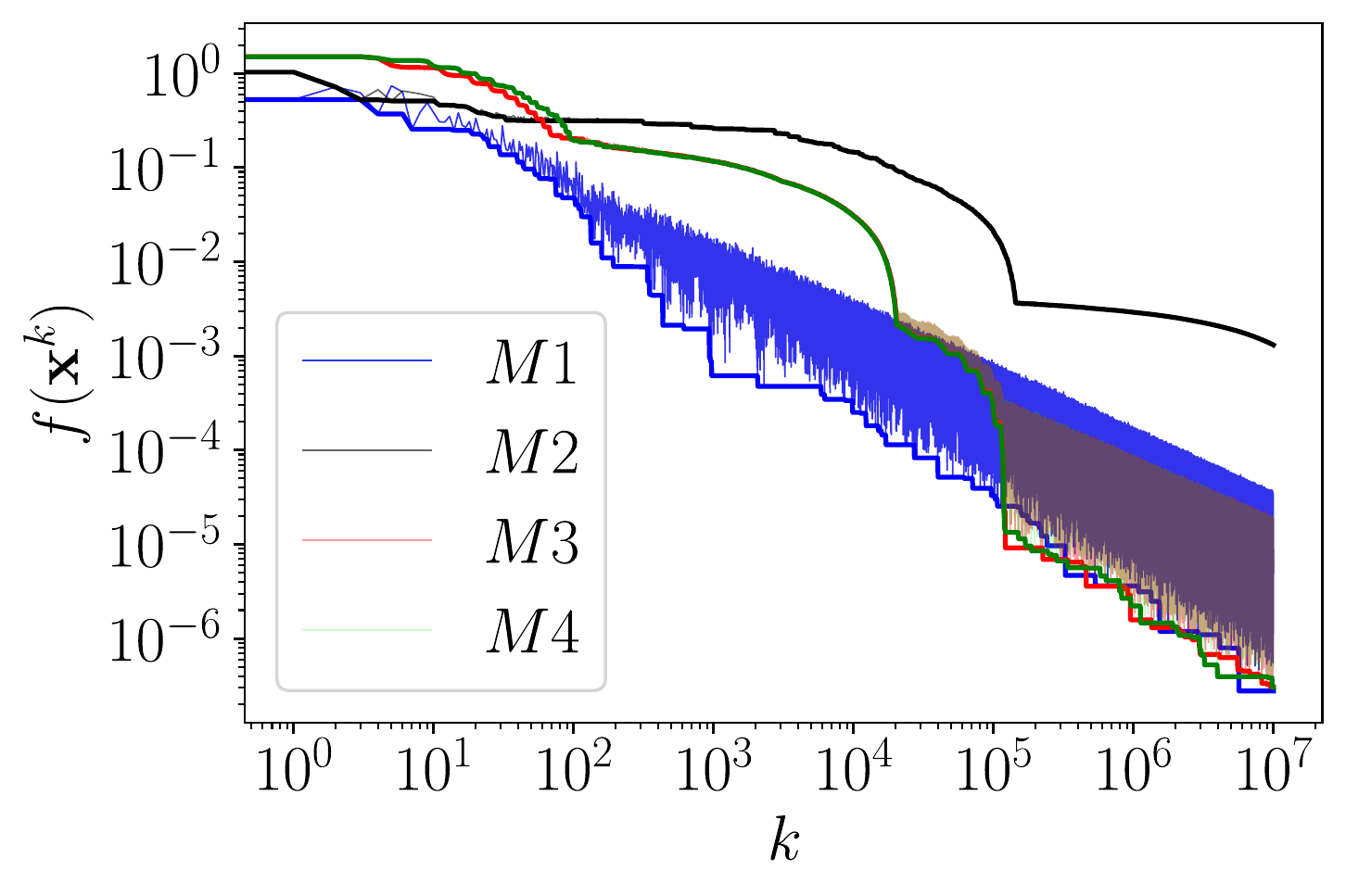} \\
(d) Test $4$ & (e) Test $5$ & (f) Test $6$
\end{tabular}
\caption{Simulation results for the tests from $1$ to $6$ using a diminishing stepsize rule are shown in (a)--(f). For each method, a thicker line that follows the graph $k \times f(\vetx^k)$ is shown, representing the smallest $f(\vetx^k)$ obtained up to iteration $k$.} \label{tests1-8}
\end{center}
\end{figure*}
\begin{figure*}[ht!]
\begin{center}
\begin{tabular}{ccc}
\includegraphics[scale=0.37]{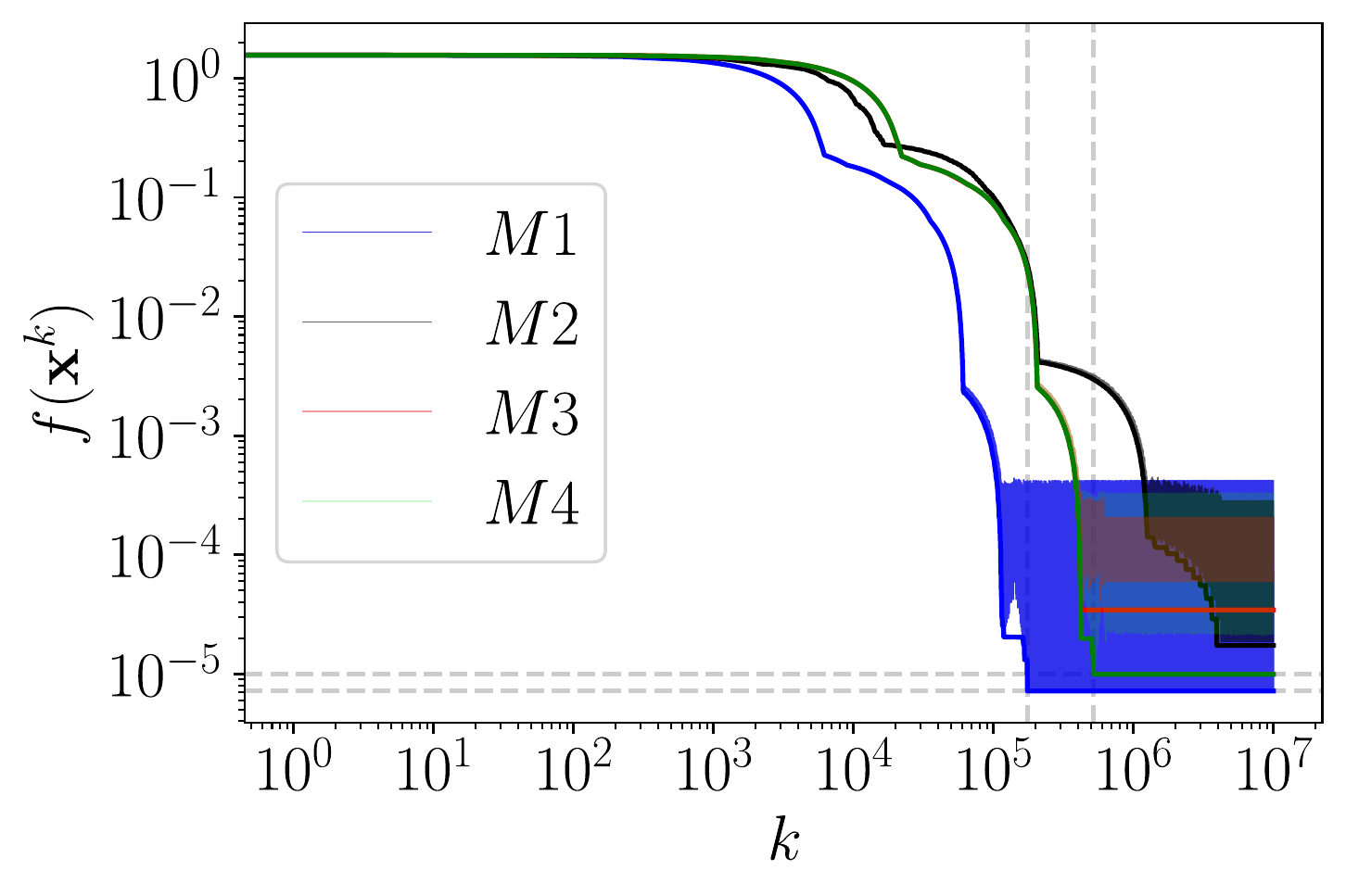} &
\includegraphics[scale=0.37]{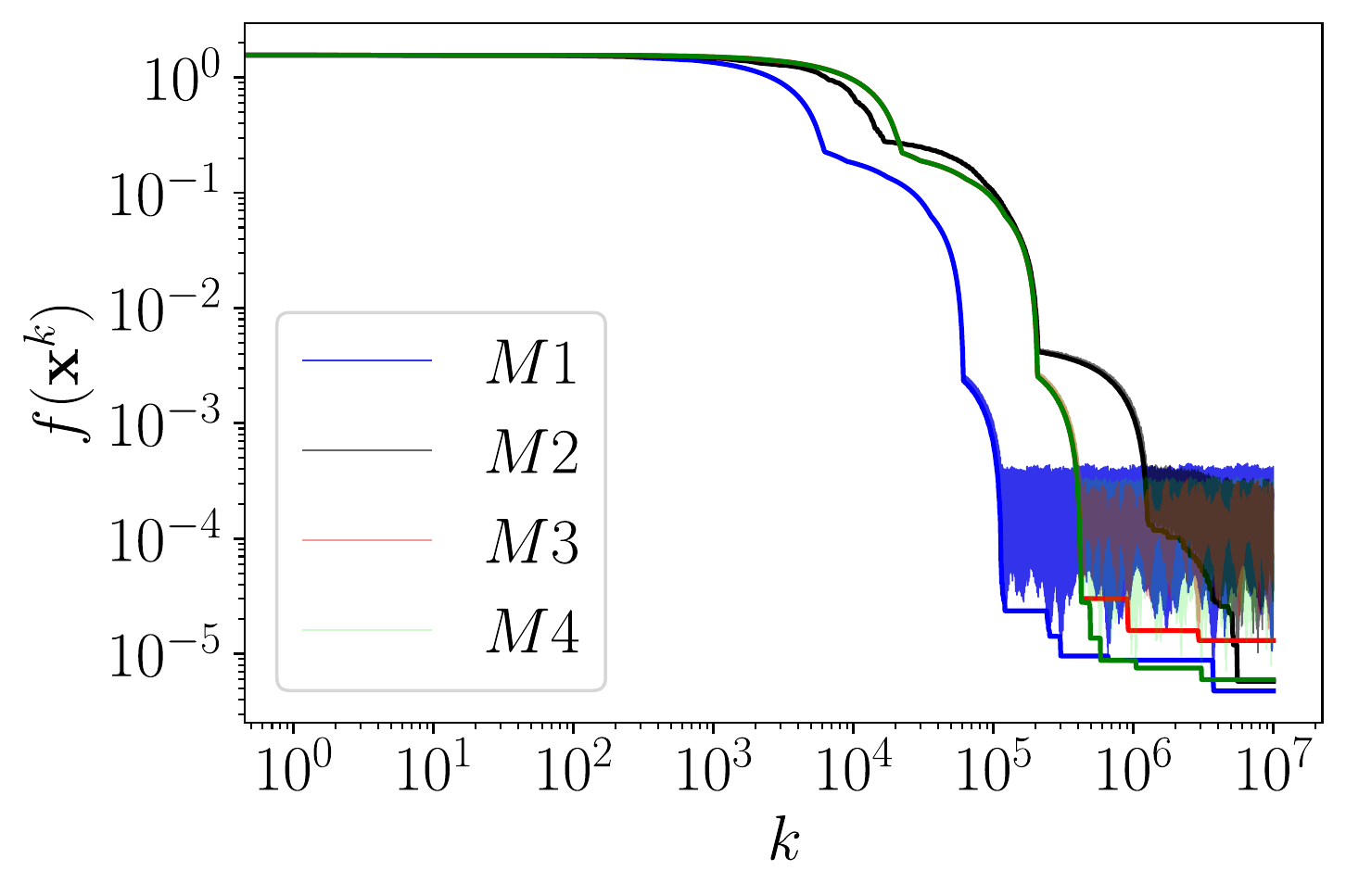} &
\includegraphics[scale=0.37]{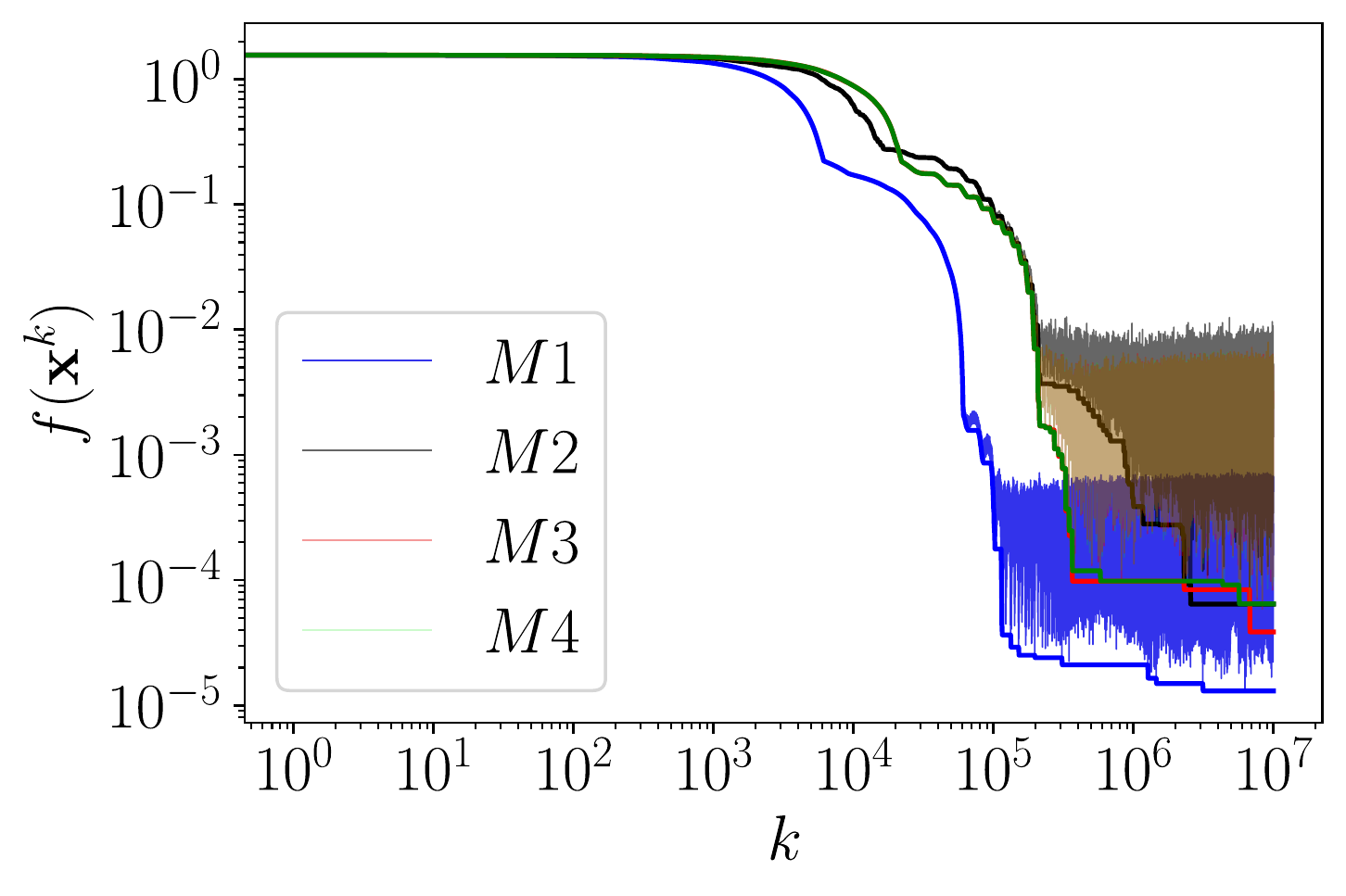} \\
(a) Test $1$ & (b) Test $2$ & (c) Test $5$
\end{tabular}
\caption{Simulation results for the tests $1$, $2$ and $5$ using $\lambda_k \equiv \lambda > 0$ are shown in (a)--(c). The values of $\lambda$ that provided the best results for each method in test $1$ (and then were also used in tests $2$ and $5$) were: $\lambda = 5 \times 10^{-4}$ for $M1$ and $\lambda = 10^{-3}$ for $M2$, $M3$ and $M4$. For each method, a thicker line that follows the graph $k \times f(\vetx^k)$ is shown, representing the smallest $f(\vetx^k)$ obtained up to iteration $k$.} \label{tests-const}
\end{center}
\end{figure*}

\rafaelb{Figure \ref{tests1-8} illustrates the behavior of $f(\vetx^k)$ for the methods $M1$, $M2$, $M3$ and $M4$ in tests from $1$ to $6$ by using diminishing stepsizes. For all methods, we draw a thicker line that follows the graph $k \times f(\vetx^k)$, representing the smallest $f(\vetx^k)$ obtained up to iteration $k$. In Fig. \ref{tests1-8}--(a), we can notice that $M1$ produces a sequence $\{\vetx^k\}$ such that $f(\vetx^k)$ decreases faster in the initial iterations in comparison to the other methods. When observing a longer time horizon, $M3$ and $M4$ produce $\{\vetx^k\}$ such that $f(\vetx^k)$ are closer to the values obtained by $M1$. Apparently, the choice of $P_{M2}$ does not produce good results for $M2$. Other neighborhood schemes could be explored, or even other transition matrices (such as min--equal neighbor scheme or weighted Metropolis--Hastings scheme that also satisfy Assumption 5 in~\cite{sundhar2009incremental}). \rafaelc{For $M1$, $M3$ and $M4$, we draw a dashed horizontal line for the first occurrence of $f(\vetx^k) < 1.0$e$-03$ (note that these lines are very close). For $M1$, this occurs 
at iteration $k = 1{,}955$ and $0.12$ seconds of CPU time. For $M3$, this occurs 
at iteration $k = 101{,}303$ and $3.34$ seconds of CPU time. 
For M4, this occurs 
at iteration $k = 98{,}833$ and $3.34$ seconds of CPU time.} Dashed vertical lines passing through such values of $k$ are also shown in the Fig. \ref{tests1-8}--(a). In addition, we insert a dashed horizontal line at the smallest value of $f(\vetx^k)$ obtained by $M4$ during all optimization process, namely $f_{M4}^{\mathrm{min}} \approx 8.7126$e$-07$, and a vertical dashed line in the iteration where this occurs, namely $k = 9{,}969{,}686$. Such solution was obtained with CPU time equal to $337.72$ seconds. A horizontal dashed line passing through the first occurrence of $f(\vetx^k) < f_{M4}^{\mathrm{min}}$ for the method $M1$ is also included. This occurs with $f(\vetx^k) \approx 8.0886$e$-07$ at iteration $k = 5{,}122{,}533$ with CPU time equal to $267.23$ seconds.}

\rafaelb{In Fig. \ref{tests1-8}--(b), a similar behavior occurs for the initial iterations: $M1$ generates a sequence $\{\vetx^k\}$ such that $f(\vetx^k)$ decreases faster in comparison to the values obtained by the other methods. However, after approximately $40{,}000$ iterations, $M3$ and $M4$ produce sequences $\{\vetx^k\}$ with smaller values for $f(\vetx^k)$ compared to the values generated by $M1$. At the end of $10{,}000{,}000$ iterations, $M3$ and $M4$ reach values of $f(\vetx^k)$ slightly smaller than those produced by $M1$.}

\rafaelb{Note that tests $1$ and $2$ are perfectly compatible with what is predicted by Theorem~\ref{teorema1} for all considered methods $M1$--$M4$. In order to explore different types of random errors, tests from $3$ to $6$ shown in Fig. \ref{tests1-8} (c)--(f), unlike tests $1$ and $2$, aim to visualize the behavior of incremental subgradient methods when Assumption~\ref{assump3} does not hold. Such situations are interesting because they can be \rafaelc{related to practical contexts where the error does not tend to zero}. Furthermore, they can be studied from the point of view of approximating the sequence $\{\vetx^k\}$ to some neighborhood of $\X^{\ast}$ when we use a constant stepsize strategy and adopt $\nu_k \equiv \nu > 0$ in Assumption~\ref{assump2}. Note that tests $3$ and $4$ have the same type of error, however the error in test $4$ is one order of magnitude smaller. In this sense, the value of $\nu$ can be chosen smaller in test $4$ and, by Theorem~\ref{teorema2}, a better error bound can be obtained. Despite the difficulties imposed on tests $3$ and $4$, which clearly impact the performance of the methods running with diminishing stepsizes, $M1$ achieved an interesting performance being able to provide more accentuated decrease in the objective function compared to the other methods. In both cases, $M1$ generates a sequence $\{\vetx^k\}$ such that $f(\vetx^k)$ is decreased up to a certain number of iterations, after which the sequence $\{\vetx^k\}$ produces less oscillation for $f(\vetx^k)$ and stabilizes. Such evidence suggests that MISSA benefited from the action of the Markov chains used, allowing indexes of rows of $A$ with higher norm (in the  subgradient calculation) to be chosen with higher frequency over the iterations, and this might have accelerated the method.}

\rafaelb{In tests $5$ and $6$ shown in Fig. \ref{tests1-8} (e)--(f), all methods performed better than in tests $3$ and $4$. In test $5$, the curve with the smallest values obtained for $f(\vetx^k)$ with the sequence $\{\vetx^k\}$ generated by $M1$ remains below the other curves generated by $M2$, $M3$ and $M4$ during the entire execution. In test $6$, similarly to what occurred in test $2$, $M1$ has better performance in the initial iterations. After about $10{,}000$ iterations, methods $M3$ and $M4$ generate $\{\vetx^k\}$ with $f(\vetx^k)$ smaller compared to the values obtained by $M1$. However, until the end of the execution, all three methods achieve similar values for $f(\vetx^k)$.}


\rafaelb{Figure~\ref{tests-const} shows the performance of the methods during $10{,}000{,}000$ iterations using a constant stepsize $\lambda > 0$ in three of the tests described in Table~\ref{table1}, namely tests $1$, $2$ and $5$. We choose $\lambda \in \{10^{-4}, 5 \times 10^{-4}, 10^{-3}, 5 \times 10^{-3}, 10^{-2}, 5 \times 10^{2}, 10^{-1}\}$ for each method, following the criterion of using the one that generates the best decrease for $f(\vetx^k)$ until the end of the $10{,}000{,}000$ iterations in test $1$. In all tests, $M1$ seems to approach a neighborhood of $\X^{\ast}$ faster compared to the other methods. In test $1$, we compare the minimum value reached by $f(\vetx^k)$ for the methods $M1$ and $M4$ and a dashed vertical line in the iteration $k$ where such a value occurs. For $M1$, this occurs with $f(\vetx^k) \approx 7.1875$e$-06$ at iteration $k = 174{,}492$ and $8.92$ seconds of CPU time. For $M4$, this occurs with $f(\vetx^k) \approx 9.9261$e$-06$ at iteration $k = 515{,}357$ and $17.71$ seconds of CPU time.}

\section{Conclusions} \label{sec.5}

In this paper, we have introduced a new method for minimizing a weighted sum of convex functions based on incremental subgradient algorithms with subgradients chosen through Markov chains.
\eduardo{No special condition is imposed on the Markov chain, allowing for transient states and
periodicity, which makes MISSA flexible enough as to contain
other methods, as the incremental cyclic subgradient method and the other methods
considered in Section IV.
This paves the way for exploring new situations
that might benefit incremental methods.}
Certain large-scale optimization problems can be addressed with the method we have presented, especially taking advantage of the flexibility of parallel processing.
Future work may move toward exploring the inclusion of transient states in the construction of the transition probability matrix and using them over a given time horizon (which is established through their own transition probabilities) in order to accelerate the convergence without changing the
objective function $f$. \eliasb{Also interesting is to investigate how theoretical convergence rates can be optimized and how to build more efficient Markov chains for incremental methods.}

\section*{Acknowledgment}
This work was supported by CNPq Grant 310877/2017-2, CNPq-Universal 421486/2016-3, FAPESP 2017/20934-9 and FAPESP-CEPID 2013/07375-0.

\appendix

\section{Preliminary results}

In this appendix we collect preliminary claims used in the convergence analysis and, when necessary, their proofs.

\begin{lemma} \emph{(Robbins-Siegmund, \cite{polyak87} - Lemma 11, p.50)} \label{lema2} Let $(\Omega, F, \mathbb{P})$ be a
probability space and $\left\{F_k \right\}$ a sequence of sub-$\sigma$ algebras of $F$. Let $\{a_k\}$, $\{c_k\}$ and
$\{\rho_k\}$ be non-negative random sequences and let $\{b_k\}$ be a deterministic sequence. Suppose that $\sum_{k=0}^{\infty} b_k < \infty$,
$\sum_{k=0}^{\infty} \rho_k < \infty$ and
$$ \EE \left[ a_{k+1} | F_k\right] \leq (1+ b_k) a_k - c_k + \rho_k,$$
holds with probability 1. Then, with probability 1, the sequence $\left\{ a_k \right\}$ converges to a non-negative
random variable and $\sum_{k=0}^{\infty} c_k < \infty$.
\end{lemma}


\begin{proposition} \label{prop1}
Let $A \in \mathbb{R}^{m \times m}$ and consider the linear system $\vetz_{k+1} = A \vetz_k$.
If for a given initial condition $\vetz_0$ the sequence $\{\vetz_k\}$ converges to zero, then there
exist $\alpha,\beta > 0$, such that
\begin{equation*}
\|\vetz_k \| \leq \alpha e^{-\beta k} \| \vetz_0 \|.
\end{equation*}
\end{proposition}
\begin{proof}
Consider the subset $E^s$ of $\mathbb{R}^m$ comprised of all initial conditions $\vetz_0$ such that the
solution $\{\vetz_k\}$ of the linear system converges to zero. $E^s$ is a vector subspace because the
system is linear, and it is $A$-invariant because the system is time-invariant so
its solution starting from any such $\vetz_k$ converges to zero.
Let us ``restrict'' $A$ into  $E^s$ by constructing a (possibly)
lower dimensional linear mapping $\tilde A: E^s\rightarrow E^s$
such that $\tilde A x = A x$ for all $x\in E^s$. The linear system $\vetz_{k+1}=\tilde A \vetz_k$, $\vetz_k\in E^s$,
constructed in this way is asymptotically stable.
From the equivalence between asymptotic and exponential stability  given in
\cite[Theorem 8.4]{hespanha2009} we have that there exist positive constants $\alpha$ and $\beta$ such that
\begin{equation*}
\| \tilde A^{k} \| \leq \alpha e^{- \beta k}, \qquad \mbox{for all} \quad k \geq 0.
\end{equation*}
This yields
\begin{equation*}
\|\vetz_k \| = \| A^k \vetz_0 \| = \| \tilde A^k \vetz_0 \| \leq \alpha e^{-\beta k} \| \vetz_0 \|,
\end{equation*}
for all $k \geq 0$.
\end{proof}

\begin{lemma} \label{lema3}
Consider an aperiodic finite state space Markov chain and let $\boldsymbol \pi^{\infty}(\boldsymbol \pi^0)=\lim_{k\rightarrow \infty} (\boldsymbol\pi^0)^T P^k$
be the limiting distribution for a given initial distribution $\boldsymbol \pi^0$.
Then, there exist $\bar{\alpha}, \bar{\beta} > 0$ such that
for any initial distribution $\boldsymbol \pi^0$
we have that
$\boldsymbol \pi^k$ converges exponentially to $\boldsymbol \pi^{\infty}(\boldsymbol \pi^0)$
in the sense that
\begin{equation*}
\| d_k \| \leq \bar{\alpha} e^{- \bar{\beta} k},
\end{equation*}
where  $d_k := \boldsymbol \pi^k - \boldsymbol \pi^{\infty}(\boldsymbol \pi^k)$.
\end{lemma}
\begin{proof}
The existence of $\boldsymbol \pi^{\infty}(\boldsymbol \pi^0)$ for aperiodic Markov
chains with finitely many states follow directly from the computation of
$\lim_{k\rightarrow\infty} P^k $ via Theorem 5.3.2 and Corollary 2.11 in \cite{cinlar2013introduction}:
$\boldsymbol \pi^{\infty}(\boldsymbol \pi^0)= \boldsymbol \pi^0 \lim_{k\rightarrow\infty} P^k.$
We shall write $d_k(\boldsymbol \pi^0) =
\boldsymbol \pi^k(\boldsymbol \pi^0) - \boldsymbol \pi^{\infty}(\boldsymbol \pi^0)$ to emphasize the dependence on the initial distribution. Then, we have that $d_k(\boldsymbol \pi^0) \to 0$ as
$k \to \infty$, and since $d_{k+1}(\boldsymbol \pi^0)^T = d_k(\boldsymbol \pi^0)^{T}P$,
Proposition \ref{prop1} \rafael{(taking $A^T = P$)} provides
\begin{equation*}
\|d_k(\boldsymbol \pi^0) \|_{\infty} \leq \alpha e^{- \beta k} \| d_0(\boldsymbol \pi^0) \|_{\infty} \leq \alpha e^{-\beta k},
\end{equation*}
where $\alpha, \beta > 0$ may depend on $\boldsymbol \pi^0$ (notice that the result in Proposition \ref{prop1} does not depend on the choice of the norm). We want to obtain
a bound that \elias{is} independent of choice of $\boldsymbol \pi^0$. For that, let $\{\vete_1, \dots, \vete_m\}$ be the canonical basis of $\mathbb{R}^m$. Again, Proposition
\ref{prop1} ensures that there are $\alpha_i, \beta_i > 0$ such that
\begin{equation} \label{aux.lema3}
\|d_k(\vete_i)\|_{\infty} \leq \alpha_i e^{- \beta_i k}, \qquad i=1,\dots,m,
\end{equation}
where $\alpha_i, \beta_i$ are parameters linked with the initial distribution $\vete_i$. Further, since $\boldsymbol \pi^k(\boldsymbol \pi^0)^T = (\boldsymbol \pi^0)^{T} P^k$ and
$\boldsymbol \pi^{\infty}(\boldsymbol \pi^0)^T = (\boldsymbol \pi^0)^{T} \lim_{k \to \infty} P^k$, 
then we have for any $\boldsymbol \pi^0 = \pi_{1}^{0} \vete_1 + \dots +
\pi_{m}^{0} \vete_m$
\begin{align*} d_k(\boldsymbol \pi^0)^T & = (\pi_{1}^{0} \vete_1 + \dots + \pi_{m}^{0} \vete_m)^{T} P^k - (\pi_{1}^{0} \vete_1 + \dots + \pi_{m}^{0} \vete_m)^{T} \lim_{k \to \infty} P^k \\
                       & = \pi_{1}^{0} \vete_1^{T} P^k + \dots + \pi_{m}^{0} \vete_m^{T} P^k - \left(\pi_{1}^{0} \vete_1^{T} \lim_{k \to \infty} P^k + \dots + \pi_{m}^{0} \vete_m^{T} \lim_{k \to \infty} P^k \right) \\
                       & = \boldsymbol \pi^k(\pi_{1}^{0} \vete_1)^T - \boldsymbol \pi^{\infty}(\pi_{1}^{0} \vete_1)^T + \dots + \boldsymbol \pi^k(\pi_{m}^{0} \vete_m)^T - \boldsymbol \pi^{\infty}(\pi_{m}^{0} \vete_m)^T \\
                       & = \sum_{i=1}^{m} d_k(\pi_{i}^{0} \vete_i)^T.
\end{align*}
Therefore, by the previous equality and equation (\ref{aux.lema3}), we have
\begin{align*} \|d_k(\boldsymbol \pi^0)\|_{\infty} & = \Bigl\| \sum_{i=1}^{m} d_k(\pi_{i}^{0} \vete_i) \Bigr\|_{\infty} \\
& \leq \sum_{i=1}^{m} |\pi_{i}^{0}| \|d_k(\vete_i)\|_{\infty} \\
& \leq \sum_{i=1}^{m} |\pi_{i}^{0}| \alpha_i e^{-\beta_i k} \leq \sum_{i=1}^{m} |\pi_{i}^{0}| \max_{j=1,\dots,m}(\alpha_j)e^{\max\limits_{j=1,\dots,m}(-\beta_j) k} \\
& \leq m \max_{j=1,\dots,m}(\alpha_j)e^{\max\limits_{j=1,\dots,m}(-\beta_j) k}.
\end{align*}
Since $\boldsymbol \pi^0$ was taken arbitrarily, the result follows taking $\bar{\alpha} = m \max_{j=1,\dots,m}(\alpha_j)$ and $\bar{\beta} = \eliasc{\max_{j=1,\dots,m}(-\beta_j)}{\min_{j=1,\dots,m}(\beta_j)}$.
\end{proof}

\begin{corollary}\label{cor-transient}
There exist $\bar{\alpha}$ and  $\bar{\beta}$ such that

\emph{(i)} $\text{Prob}\left( s_\ell(\delta \eliasc{n(k)}{k} )\in\mathcal T \right) \leq \bar{\alpha} e^{- \bar{\beta} k}$ and

\emph{(ii)} $ \text{Prob}\left (s_{\ell}(\delta \eliasc{n(k)}{k}) \in \mathcal{R}_v  \right)
[ \boldsymbol \pi_{\delta,v}]_i - [\boldsymbol \pi_{\ell}^{\infty} ]_i \leq \bar{\alpha} e^{- \bar{\beta} k}$.
\end{corollary}
\begin{proof}
Lemma \ref{lema3} applies to $P^\delta$
because the $\delta$-step chain $\{s_\ell(\delta k), k\geq 0\}$ is aperiodic.
The limiting probability of a transient state $i$ is always zero,
that is, $[\lim_{k\rightarrow \infty} (\boldsymbol\pi^0_{\eliasc{}{\ell}})^T P^{\delta k}]_i=0$,
and Lemma \ref{lema3} establishes that convergence is exponentially fast
(with uniform parameters $\bar{\alpha}$ and  $\bar{\beta}$), which leads to (i).
Similarly, $\text{Prob}\left (s_{\ell}(\delta \eliasc{n(k)}{k}) \in \mathcal{R}_v  \right)$ converges
exponentially fast to a limiting constant
$c = \sum_{i\in \mathcal{R}_v }[\lim_{k\rightarrow \infty} (\boldsymbol\pi^0_{\eliasc{}{\ell}})^T P^{\delta k}]_i$,
and $c[ \boldsymbol \pi_{\delta,v}]_i - [\boldsymbol \pi_{\ell}^{\infty} ]_i =0$, and some algebra leads
to (ii).
\end{proof}

We point out that the Markov chain with
transition probability matrix $P^\delta$ is an aperiodic chain \cite[Theorem 3.7]{cinlar2013introduction} and
in this case $\Delta = \lim_{k\rightarrow\infty} P^{\delta k} $ can be computed
via Theorem 5.3.2 and Corollary 2.11 in \cite{cinlar2013introduction}.

\begin{lemma} \label{lema4} Consider a time-homogeneous Markov chain with finite state space and probability transition matrix $P$. Let $\delta > 0$ be the period of this chain, \rafaelc{$\Phi(k,t) = P^{k-t}$ with $k>0$, $t \geq 0$, $k > t$ and $\Delta = \lim_{k \to \infty} P^{\delta k}$}. Then, there exist $\alpha, \beta > 0$ such that
\begin{description}
 \item[(i)] $\| \Phi(\delta k, \delta t) - \Delta \|_{\infty} \leq \alpha e^{-\beta(k-t)}$;
 \item[(ii)] $\| \Phi(\delta k + i, \delta t + j) - \Delta P^{(i-j)} \|_{\infty} \leq \alpha e^{-\beta(k-t)}$ for all $i,j \geq 0$ with $i \geq j$.
\end{description}
\end{lemma}
\begin{proof}
\textbf{(i)} We denote $\boldsymbol \pi^{\delta \infty}(\boldsymbol \pi^0)^T = (\boldsymbol \pi^0)^{T} \Delta$ and thus
\begin{align*}  d_{\delta k}(\boldsymbol \pi^0)^T &= \boldsymbol \pi^{\delta k}(\boldsymbol \pi^0)^T - \boldsymbol \pi^{\delta \infty}(\boldsymbol \pi^0)^T \\
& = (\boldsymbol \pi^0)^T (P^{\delta k} - \Delta) \to 0,\end{align*}
as $k \to \infty$ and by noticing that $d_{\delta(k+1)}(\boldsymbol \pi^0)^T = d_{\delta k}(\boldsymbol \pi^0)^{T}P^{\delta}$, we can apply Proposition \ref{prop1} and Lemma
\ref{lema3} and obtain geometric convergence for the sequence $\{ d_{\delta k}(\boldsymbol \pi^0) \}$ for any $\boldsymbol \pi^0$. By observing that
$\boldsymbol \pi^{\delta k}(\vete_i)$ is the $i$-th row of the matrix $P^{\delta k}$ and $\boldsymbol \pi^{\delta \infty}(\vete_i)$ is the $i$-th row of $\Delta$, then
$d_{\delta k}(\vete_i)$ is the $i$-th row of the matrix $P^{\delta k} - \Delta = \Phi(\delta k + 1, 1) - \Delta$. Bringing this and Lemma \ref{lema3} together yields
\begin{align*} \| \Phi(\delta k + 1, 1) - \Delta \|_{\infty} & = \max_{i=1,\dots,m} \| d_{\delta k}(\vete_i) \|_1 \\
& \leq \max_{i=1,\dots,m} (m \| d_{\delta k}(\vete_i) \|_{\infty}) \\
& \leq \max_{i=1,\dots,m} (m \bar{\alpha} e^{-\bar{\beta} \delta k}) = \alpha e^{-\beta k},
\end{align*}
where $\alpha = m \bar{\alpha} > 0$ and $\beta = \delta \bar{\beta} > 0$. The result follows by noticing that $\Phi(\delta k, \delta t) = \Phi(\delta(k-t) + 1, 1)$.

\textbf{(ii)} It is sufficient note that
\begin{align*} \| \Phi(\delta k + i, \delta t + j) - \Delta P^{(i-j)} \|_{\infty} & = \| P^{\delta(k-t) + i - j} - \Delta P^{(i-j)} \|_{\infty} \\
&= \| (P^{\delta(k-s)} - \Delta) P^{(i-j)} \|_{\infty} \\
& \leq \| \Phi(\delta k, \delta t) - \Delta \|_{\infty} \| P^{(i-j)} \|_{\infty} \\
& \leq \alpha e^{-\beta(k-t)},
\end{align*}
where we use \textbf{(i)} and $\| P^k \|_{\infty} \leq 1$ for all $k \geq 0$.
\end{proof}

\bibliographystyle{plain}
\bibliography{references}

\end{document}